\documentclass[reqno]{amsart}

\usepackage[draft]{changes}
\usepackage{graphicx}
\usepackage{amsmath,amssymb}
\usepackage{amsthm}
\usepackage{xcolor}
\usepackage{mathrsfs}
\usepackage{paralist}
\usepackage{multirow}
\usepackage{mathtools}
\usepackage[T1]{fontenc}
\usepackage[utf8]{inputenc}
\usepackage[english]{babel}
\DeclareMathAlphabet\mathbfcal{OMS}{cmsy}{b}{n}
\usepackage[figurename=Figure]{caption}
\usepackage{lineno}
\usepackage{bbm}
\usepackage[normalem]{ulem}
\usepackage[percent]{overpic}
\usepackage{color}
\usepackage{tikz}
\usepackage{yhmath}

\usepackage{wrapfig}
\usepackage{float}
\usepackage{booktabs}
\usepackage{makecell}
\usepackage{subcaption}
\usepackage{comment}
\usepackage{cancel}
\theoremstyle{definition}
\newtheorem{definition}{Definition}[section]

\usepackage{hyperref}
\hypersetup{
    colorlinks=true,     
    linkcolor=red,       
    citecolor=green,     
    filecolor=magenta,   
    urlcolor=blue        
}

\definechangesauthor[name={Luca}, color=orange]{LM}

\newcommand{\eg}{\emph{e.g.}} 

\newcommand{\del}[1]{}

\usepackage{authblk}

\usepackage{amsaddr}
\usepackage{siunitx}

\usepackage[a4paper,bindingoffset=0.0in,%
            left=1.0in,right=1.0in,top=1.in,bottom=1.in,%
            footskip=.4in]{geometry}

\theoremstyle{remark}

\title[Symmetric CAEs: enhancing latent stability via differential geometry]{Convolutional Symmetric AutoEncoders: enhancing latent stability via differential geometry}

\author{Gaspare Li Causi$^{1}$, Niccolò Tonicello$^{1}$, Luca Magri$^{2,3}$ and Gianluigi Rozza$^{1}$}
\address{$^1$ mathLab, Mathematics Area, SISSA, via Bonomea 265, I-34136 Trieste, Italy}
\address{$^2$ Department of Aeronautics, Imperial College London, SW7 2AZ, London, UK}
\address{$^3$ Politecnico di Torino, DIMEAS, Corso Duca degli Abruzzi, 24 10129 Torino, Italy}

\begin{document}
\begin{abstract}
Autoencoders (AEs) have emerged as powerful tools for non-linear dimensionality reduction, often surpassing traditional linear methods such as Proper Orthogonal Decomposition (POD) in scenarios characterized by slowly decaying Kolmogorov $n$-widths. In the realm of Reduced-Order Modelling (ROM), these models are increasingly utilized to learn low-dimensional representations of solution manifolds associated with parametric Partial Differential Equations (PDEs). However, the high expressivity of AEs presents a challenge: although trained networks typically minimize reconstruction error, they often struggle to capture the essential properties necessary for building accurate and robust ROMs.
Recent works by \cite{otto2023learning} and \cite{brivio2025deep} have tackled this challenge in fully connected AEs by proposing representation-consistent architectures, which preserve some of the properties belonging to POD. This study builds upon that concept by extending representation consistency for convolutional layers. We introduce a novel class of symmetric Convolutional AutoEncoders (CAEs) designed to embody the primary properties of manifold parametrization mappings. When integrated into a ROM framework, this architecture demonstrates significantly improved predictive capabilities. Specifically, we compared the performance of the ROMs based on classical and symmetric CAEs on three one dimensional academic test cases, namely the Linear Advection, Viscous Burger and Kuramoto–Sivashinsky equation. Numerical results demonstrate that our proposed symmetric approach consistently yields more accurate latent trajectories, lower reconstruction errors, and enhanced model robustness.  
\end{abstract}

\maketitle

\section{Introduction}

Over the past few decades, numerical simulations of physical systems have provided invaluable insight, significantly enhancing our understanding of complex phenomena. These systems are typically modeled using Partial Differential Equations (PDEs) for which closed-form solutions are rarely attainable. Instead, they must be solved through high-fidelity numerical frameworks such as Finite Difference (FD), Finite Element Methods (FEM), and Finite Volume (FV) techniques.

The complexity of these systems does not translate directly into a more difficult implementation of those numerical schemes, but rather arises from the challenge of managing diverse temporal and spatial scales. These scales, often differing by several orders of magnitude, do not evolve separately; they interact dynamically due to the couplings arising from the nonlinear terms, leading to phenomena such as the energy cascade observed in viscous dissipative systems like turbulence. To achieve accurate solutions, a simulation must resolve the entire spectrum of scales, from the largest structures down to the smallest dissipative eddies. This requirement results in discretized systems with a massive number of Degrees of Freedom (DoF).

Despite the continuous advancement of High-Performance Computing (HPC) facilities over the last thirty years \cite{hey1997high, raj2020high}, certain realistic solutions remain virtually unattainable within reasonable time frames. A prominent example is the Direct Numerical Simulation (DNS) of an entire aircraft under cruise conditions, which remains computationally prohibitive even with modern HPC hardware. Even in scenarios where DNS is technically feasible, its application is often impractical for multi-query contexts such as structural parameter optimization or for real-time prediction models embedded in control systems.

In this context  Reduced Order Models (ROMs) comes into play. Their objective is to speed up the numerical solutions, hence drastically decreasing the computational burden, at the cost of losing some degree of accuracy \cite{quarteroni2014reduced}. 
Within this framework, classical turbulence modeling approaches such as Large Eddy Simulation (LES) \cite{sagaut2006large} and Reynolds-Averaged Navier-Stokes (RANS) \cite{wilcox1998turbulence} can be interpreted as specialized forms of model reduction. In these approaches, only the most energetic large-scale structures of the flow are explicitly resolved, while the influence of the unresolved small scales is accounted for through various closure models. This simplification enables the use of coarser computational grids, significantly reducing the required degrees of freedom and, consequently, the time-to-solution.

Dissipative systems, which arise from nonlinear PDEs, exhibit a notable property: their long-time dynamics reside on a simpler, finite-dimensional structure known as the \textit{attractor} \cite{temam2012infinite}. The \textit{attractor} is an invariant subset with respect to the semigroup described by the solution map of the dissipative systems, thus encapsulating all essential information regarding the non-trivial dynamics of the system \cite{temam1990inertial}. Given that this \textit{attractor} is finite-dimensional, the long-term dynamics following transient behaviors can be effectively represented by a system of Ordinary Differential Equations (ODEs) \cite{zelik2014inertial}.
This finite-dimensional aspect can be articulated through the concept of an \textit{inertial manifold}, defined as a smooth finite-dimensional manifold $\mathcal{M}$ that remains invariant under the semigroup associated with the solution operator. Additionally, the inertial manifold contains the global attractor and attracts all trajectories at an exponential rate \cite{foias1988inertial, temam1990inertial, zelik2014inertial}. 
For certain dissipative dynamical systems, such as the Navier-Stokes equations, reaction-diffusion equations, and damped hyperbolic equations, the existence of an \textit{inertial manifold} has been established \cite{babin1992attractors, chepyzhov2002attractors}. The presence of these inertial manifolds thus lends rigorous theoretical support to the \textit{Manifold Hypothesis}, which posits that the set of all solutions of a parametric PDE lies on, or can be well-approximated by, a low-dimensional smooth manifold $\mathcal{M}$ embedded within the high-dimensional state space \cite{whiteley2025statisticalexplorationmanifoldhypothesis, temam1990inertial}. Within this framework, Reduced Order Models (ROMs) seek to identify a mapping from the high-dimensional physical space, where the manifold is embedded, to a low-dimensional representation \cite{quarteroni2014reduced} commonly referred to in ROM terminology as the \textit{latent space} \cite{racca2023predicting, fresca2021comprehensive, brivio2025deep}.

Once a proper reduced representation is established, a latent dynamical model is required to evolve the low-dimensional state over time \cite{farenga2025latent}. The low dimensionality of the latent space is what makes the time integration of the reduced system computationally cheaper, allowing for a much faster solution of the physical problem. Traditionally, Reduced Basis (RB) methods approximate this solution manifold as a linear combination of basis functions, typically derived via Proper Orthogonal Decomposition (POD), Dynamic Mode Decomposition (DMD) \cite{schmid2010dynamic}, or Greedy sampling algorithms.

However, linear RB methods face significant limitations, most notably the slow decay of the Kolmogorov $n$-width for advection-dominated problems \cite{arbes2025kolmogorov}. In such cases, the linear subspace spanned by the reduced basis fails to efficiently represent the nonlinear manifold, requiring an excessively large number of basis functions to maintain acceptable accuracy and thereby diminishing the computational advantage of the reduction.


The advent of Machine Learning has significantly transformed the landscape of Model Order Reduction (MOR), providing powerful tools to circumvent the restrictions inherent in traditional Reduced Basis methods \cite{peherstorfer2022breaking}. By utilizing specialized Neural Network architectures, one can establish nonlinear mappings that connect the high-fidelity full-order space with a low-dimensional latent representation. These nonlinear mappings facilitate superior information compression compared to linear frameworks, thereby addressing the challenges posed by the slow decay of the Kolmogorov $n$-width in advection-dominated regimes \cite{peherstorfer2022breaking}. Currently, several architectures are prominently used for this purpose, most notably Graph Neural Networks (GNNs) \cite{pichi2024graph, barwey2023multiscale}, fully connected networks (MLPs) \cite{brivio2025deep}, Convolutional Neural Networks (CNNs) \cite{racca2023predicting, ozalp2025stability, maulik2021reduced, romor2023non}, and Variational Autoencoders (VAEs) \cite{tonioni2026vivaldy, solera2024beta, zighed2025uncertainty}.


Regarding the temporal evolution of the system, a wide range of neural architectures have been explored to approximate the latent dynamics. These range from traditional Recurrent Neural Networks (RNNs), such as Long Short-Term Memory (LSTM) and Gated Recurrent Units (GRU) \cite{vlachas2018data, gonzalez2018deep, fresca2023long}, to more recent frameworks like Neural Ordinary Differential Equations (Neural ODEs) \cite{chen2018neural, dutta2021neural}, Koopman-based operators \cite{gupta2025mori, menier2025interpretable} and Transformers \cite{tonioni2026vivaldy, solera2024beta}. For a comprehensive overview of these methodologies and additional references, the reader is referred to \cite{farenga2025latent}. 


The governing equations of the latent dynamics can be rigorously derived under the assumption that the solution set of the parametric PDE constitutes a differentiable manifold. This derivation leverages the properties of manifold parametrization via diffeomorphic mapping which, when restricted to the manifold, are differentiable, invertible, and have a differentiable inverse too. In this framework, if $\mathbf{u} \in \mathcal{M}$ represents the high-fidelity state on the manifold $\mathcal{M}$, and $\boldsymbol{z} \in U_\alpha$ is its latent representation, the latent velocity $\dot{\boldsymbol{z} \in U_\alpha}$ is formally linked to the physical velocity $\mathbf{f}(\mathbf{u})$ through the Jacobian of the mapping.

While this framework is mathematically solid, a significant challenge arises in practice: in standard autoencoder architectures, the decoder $D$ does not strictly serve as the exact inverse of the encoder $E$ (i.e., $E \circ D \neq I$). This discrepancy introduces a fundamental structural inconsistency that can invalidate the mathematically derived relations governing the latent evolution. It is important to highlight that this issue, intrinsic to the approximation of non-linear mapping through neural networks, is absent in linear RB-POD methods, where the projection and reconstruction operators are duals by construction. 

Although recent literature has proposed fully connected AEs designed to enforce the invertibility required for consistent manifold parametrization  \cite{brivio2025deep, otto2023learning}, this kind of architecture faces significant practical hurdles. Specifically, they suffer from poor scalability in terms of the number of trainable parameters relative to the input vector dimension and exhibit low sensitivity to local topological features of the physical domain. These limitations motivate the transition towards convolutional architectures \cite{goodfellow2016deep}.

In this work, we extend the symmetry properties developed by \cite{otto2023learning, brivio2025deep} to hybrid 1D convolutional autoencoders. The term hybrid refers to the structure of the network: while the encoder consists of traditional convolutional layers, similar to those previously utilized in convolutional AEs, the decoder is constructed using a fully connected layer with matrices parameterized by a user-defined number of tunable parameters. These hybrid convolutional AEs, combined with a simple and straightforward latent dynamical model, are employed to develop ROMs for 1D dynamical systems, which serve as classic academic test cases.

The remainder of this work is organized as follows. Section \ref{sec: Methodology} establishes the theoretical foundation through a formal derivation of the latent governing equations. Building upon these principles, Section \ref{sec: AE} details the construction of symmetric Convolutional Autoencoders designed to satisfy the representation consistency requirement. The integration of these architectures into a coherent Reduced Order Model is addressed in Section \ref{sec: Data-Driven}, which outlines the interaction between the AE and latent dynamical models alongside the decoupled training protocols. The practical utility of the proposed framework is evaluated in Section \ref{sec: Numerical experiments}, which describes the numerical experiments and physical test cases selected for validation. Section \ref{sec: Results} presents a comprehensive performance analysis, comparing the symmetric ROM architecture against standard model order reduction techniques. Finally, Section \ref{sec: Conclusion} summarizes the study's findings and discusses current limitations while proposing directions for future extensions.

\section{Methodology}
\label{sec: Methodology}


\subsection{Latent dynamics}

\label{sec: Latent Equation}
In this section we derive the latent governing equation as a projection of the FOM equation system, by exploiting the basics of differential geometry. To do so, we adopt the mathematical framework and notation introduced by \cite{farenga2025latent}, which is specifically tailored for applications in nonlinear model order reduction.

Consider the system of ODEs resulting from the semi-discretization of a parametric PDE, where the parameters are collected in a vector \(\boldsymbol{\mu} \in \mathcal{P} \):
\begin{equation}
\label{eq: FOM equation}
\left\{\begin{array}{l}
\dot{\mathbf{u}}_h(t ; \boldsymbol{\mu})=\mathbf{f}_h\left(t, \mathbf{u}_h(t ; \boldsymbol{\mu}) ; \boldsymbol{\mu}\right), \quad t \in\left(t_0, T\right], \\
\mathbf{u}_h\left(t_0 ; \boldsymbol{\mu}\right)=\mathbf{u}_{0, h}(\boldsymbol{\mu}).
\end{array}\right.    
\end{equation}
Here, the subscript \(h\) denotes a high-fidelity discretization method (e.g., FEM, FVM, or FD). The set of time-dependent, parametric solutions \(\mathbf{u}_h \in \mathbb{R}^{N_h}\) forms a solution manifold:
\begin{equation*}
    S_h = \left\{ \mathbf{u}_h(t; \boldsymbol{\mu}) : t \in [t_0, T], \, \boldsymbol{\mu} \in \mathcal{P} \right\} \subset \mathbb{R}^{N_h}
\end{equation*}
with \(N_h\) denoting the number of degrees of freedom of the discretized dynamical system. Assuming \(S_h\) to be a differentiable \(m\)-dimensional manifold with \(m = \dim(\mathcal{P}) + 1\) there exists a family of injective parameterizations
\begin{equation*}
\label{eq: exact decoder}
\mathbf{x}_\alpha : U_\alpha \subset \mathbb{R}^m \to S_h \subset \mathbb{R}^{N_h},
\end{equation*}
where each \(\mathbf{x}_\alpha\) maps a local coordinate \(\boldsymbol{z} \in U_\alpha\) to a point on the manifold \(S_h\). 
To provide a wider view, the map $\mathbf{x}_\alpha$ could be approximated by the decoder $D$ of an AE, while $\mathbf{x}_\alpha ^{-1}$ is approximated by the encoder $E$.

Let $\beta : I \to \mathbf{x}_\alpha (U_\alpha)$ denote the trajectory traced on the manifold by the high-fidelity solution $\mathbf{u}_h(t; \boldsymbol{\mu}^*)$ for a fixed parameter instance $\boldsymbol{\mu}^*$. This trajectory is defined on an open interval $I \subseteq (t_0, T]$ such that $\beta(t) = \mathbf{u}_h(t; \boldsymbol{\mu}^*) \in \mathbf{x}_\alpha (U_\alpha)$ for all $t \in I$. The tangent vector to $\beta$ is the temporal derivative $\dot{\beta}(t) = \dot{\mathbf{u}}_h(t; \boldsymbol{\mu}^*)$. 
If $\gamma: I \to U_\alpha$ represents the preimage of $\beta$ under the map $\mathbf{x}_\alpha$, then $\dot{\boldsymbol{z} }(t; \boldsymbol{\mu}^*)$ is the tangent vector to the curve $\gamma$ at time $t$. The structural relationship between $\gamma$, $\beta$, and $\mathbf{x}_\alpha$ is given by the composition $\beta = \mathbf{x}_\alpha \circ \gamma$. 
We denote by $d\mathbf{x}_\alpha: T_{\boldsymbol{z}} U_\alpha \to T_{\mathbf{u}_h} \mathcal{M}$ the differential of the manifold parametrization. By the definition of the differential and the application of the chain rule, we obtain:
\begin{equation}
    \dot{\mathbf{u}}_h(t; \boldsymbol{\mu}^*) = \dot{\beta}(t) = d\mathbf{x}_\alpha (\gamma(t))[\dot{\gamma}(t)] =d\mathbf{x}_\alpha  (\mathbf{z}(t; \boldsymbol{\mu}^*))[\dot{\mathbf{z}}(t; \boldsymbol{\mu}^*)].
\end{equation}
While this relation was derived for a specific trajectory associated with $\boldsymbol{\mu}^*$ , the result generalizes to any parameter $\boldsymbol{\mu}\in\mathcal{P}$.
Applying the change of variables to the FOM in Equation \ref{eq: FOM equation} yields an implicit expression for the latent time derivative:
\begin{equation*}
    d\mathbf{x}_\alpha  (\mathbf{z}(t; \boldsymbol{\mu}))[\dot{\mathbf{z}}(t; \boldsymbol{\mu})] = \mathbf{f}_h\left(t, \mathbf{x}_\alpha (\boldsymbol{z}(t ; \boldsymbol{\mu})) ; \boldsymbol{\mu}\right).
\end{equation*}
Deriving an explicit evolution law requires pre-multiplying both sides by the inverse of the differential of the manifold parametrization $(d \mathbf{x}_\alpha)^{-1}$. Diffeomorphism properties ensure that this inverse exists and is equivalent to the differential of the inverse of the manifold parametrization \cite{lee2003smooth}, namely $(d \mathbf{x}_\alpha)^{-1} = d(\mathbf{x}_\alpha^{-1})$. Projecting the high-fidelity governing equations into the latent space via such identity defines the evolution of the reduced state variables:
\begin{equation}
\label{eq: latent equation}
\dot{\boldsymbol{z}}(\boldsymbol{z},t ; \boldsymbol{\mu})=d \mathbf{x} _\alpha ^{-1} (\boldsymbol{z}(t ; \boldsymbol{\mu})) [\mathbf{f}_h\left(t, \mathbf{x}_\alpha (\boldsymbol{z}(t ; \boldsymbol{\mu})) ; \boldsymbol{\mu}\right)], \quad t \in I.
\end{equation}
The system described by Equation \ref{eq: latent equation} evolves autonomously once the map $\mathbf{x}_\alpha$ is established. Despite the reduced dimensionality of the latent state $\boldsymbol{z}$, the integration process still necessitates the evaluation of high-dimensional vectors, including the reconstructed state $\mathbf{x}_\alpha(\boldsymbol{z})$ and the full-order operator $\mathbf{f}_h$. In practice, this requirement can make the direct integration of the latent equations more computationally demanding than the original high-fidelity system \cite{romor2023non}. 

The ROM circumvents this bottleneck by approximating the right-hand side of the latent equations with a suitable NN (later named Latent Flow Network). While $\mathbf{x}_\alpha$ remains a linear map and reduces to the standard POD basis when the manifold $S_h$ is an affine space, non-linear manifolds require more sophisticated approximations. For such tasks, autoencoders have proven to be an ideal candidate for capturing complex, non-linear manifold structures.  

\section{AutoEncoders}
\label{sec: AE}
AutoEncoders constitute a broad class of neural network architectures that learn, in an unsupervised manner, to compress data into a lower-dimensional representation. The encoder, denoted by $E$, maps the input data to a latent space by discovering an appropriate reduced representation. The complementary component, the decoder $D$, approximates the inverse of this mapping and reconstructs the input from its compressed latent code, ideally recovering the original high-dimensional data. Common variants of this framework include variational Autoencoders \cite{kingma2013auto}, denoising Autoencoders \cite{vincent2008extracting}, and contractive Autoencoders \cite{rifai2011contractive}. Architecturally, these models are distinguished primarily by the types of layers they employ: convolutional AutoEncoders are built from convolutional layers, whereas fully connected AutoEncoders are composed of dense (fully connected) layers.

Owing to their high expressive power, AEs can be employed to tackle a broad range of tasks. In many applications, however, it is desirable to bias the AE toward solutions of a particular type, even when alternative representations achieve the same reconstruction error. One simple example is the reduction of an incompressible fluid velocity field, for which it is important not only that the reconstructed velocity looks similar to the original field, but also that the divergence free condition is preserved. Such preferences can be weakly enforced by modifying the loss function in an ad hoc manner \cite{olmo2022physics} or strongly by carefully constraining the architecture of the AE \cite{mohan2023embedding}. Otto et al.~\cite{otto2023learning} were the first to introduce AEs that are intrinsically representation-consistent, in the sense that any latent point remains invariant under the composition of decoder and encoder. By imposing a specific structural constraint on the AE, they effectively realized a nonlinear projector onto the slow manifold along which the reduced dynamics evolve. A subsequent theoretical analysis of these architectures was provided by Brivio et al.~\cite{brivio2025deep}, who derived upper and lower bounds on the reconstruction error in terms of the Lipschitz constants of the activation functions and POD error bounds at each encoder/decoder layer stage. Both studies \cite{brivio2025deep,otto2023learning} were restricted to fully connected AE architectures.

Let us consider an AE composed by an encoder \(E\) and a decoder \(D\). We assume that the encoder and decoder comprise the same number of convolutional layers, and we denote this common depth by \(n_d\). 
In order to preserve symmetry (\(E \circ D = id\)) the following conditions are needed \cite{brivio2025deep}:
\begin{itemize}
    \item Bi-lipschitz activation function.
    \item Orthogonal / Bi-orthogonal linear layer pairs.
\end{itemize}
\begin{figure}[h!]
    \centering
    \includegraphics[width=1\linewidth]{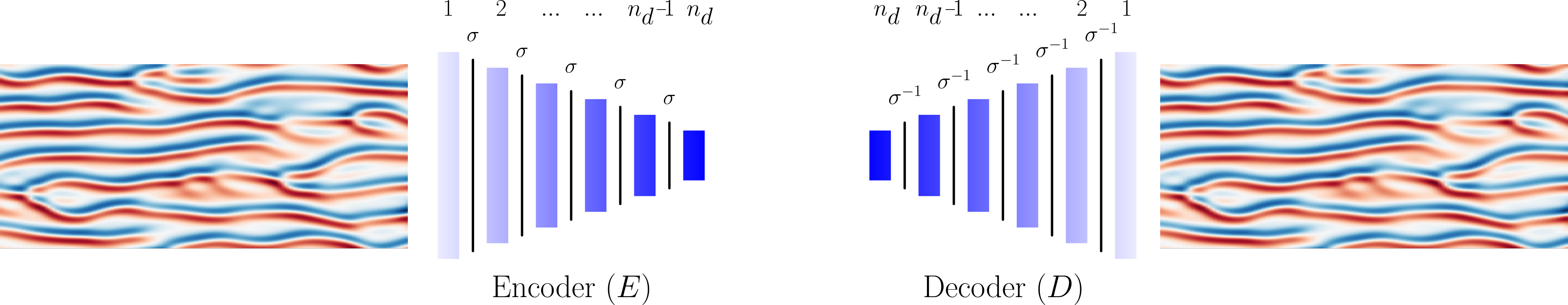}
    \caption{The architecture consists of a nonlinear encoder ($E$) for dimensionality reduction and a symmetric decoder ($D$) for physical field reconstruction. Here, $\sigma$ denotes a Bi-lipschitz nonlinear activation function applied element-wise between subsequent encoder layers, while $\sigma^{-1}$ represents its inverse, utilized within the decoder layers to maintain the structural symmetry and bi-orthogonality of the mapping.} 
    \label{fig: AE structure}
\end{figure}
Let layer $1$ denote the first layer of the encoder and layer $n_d$ denote the last layer. The opposite indexing applies to the decoder, as schematically shown in Figure \ref{fig: AE structure}. Since the decoder is a mirrored ``copy'' of the encoder, a relation must hold between the last encoder layer and the last decoder layer (according to the enumeration given in Figure \ref{fig: AE structure}), the first encoder layer and the first decoder layer, and, more generally, between the $i$-th encoder layer and the $i$-th decoder layer.
If the layers are fully connected and represented by the weight matrices  
$$
E_i \in \mathbb{R}^{n_{i+1} \times n_i}, \quad D_{i} \in \mathbb{R}^{n_i \times n_{i+1}},
$$  
with $n_{i+1} \le n_i$, then bi-orthogonality means
\begin{equation}
    E_i D_{i} = \mathbb{I} \in \mathbb{R}^{n_{i+1} \times n_{i+1}}.
    \label{eq: bi-ortho condition}
\end{equation}

Several strategies have been proposed to enforce the condition defined in Equation \ref{eq: bi-ortho condition}. For instance, Brivio and Franco~\cite{brivio2025deep} implemented a matrix decomposition approach where the encoder and decoder weight matrices are expressed as a product of five distinct matrices belonging to specific classes, such as diagonal, orthogonal, and general matrices. Alternatively, Otto et al.~\cite{otto2023learning} pursued a different methodology by defining the decoder weight matrices as the generalized inverse of the corresponding encoder matrices. It is important to note that while the Moore-Penrose pseudoinverse is unique, the generalized inverse is not, providing additional degrees of freedom in the architectural design.

If $E_i$ is an orthogonal matrix, the orthogonality condition implies that $D_{i} = E_i^T $, which automatically satisfies Equation \ref{eq: bi-ortho condition}.

While the theoretical significance of symmetric fully connected autoencoders is considerable, their practical utility is often constrained by inherent limitations. In contrast, convolutional neural networks and autoencoders have emerged as the standard for large-scale applications. With the goal of creating symmetric AEs that can make a meaningful impact in practical applications, we are taking a step forward by developing 1D convolutional symmetric AEs.

\subsection{Convolutional Symmetric AutoEncoders}
In this section, we extend the notion of symmetric AEs to convolutional architectures. For the sake of simplicity, we will assume that there are no biases involved. The inclusion of bias is straightforward, as explained in previous works on symmetric AEs \cite{otto2023learning, brivio2025deep}. 

The convolution of a function against a kernel is a linear operator. Let $\boldsymbol{k} \in \mathbb{R}^{2n_k + 1}$ denote the (discrete) kernel:
$$
\boldsymbol{k} = [k_{-n_k}, k_{-n_k +1}, ..., k_0, ... ,k_{n_k -1}, k_{n_k}].
$$
In the discrete form, the convolution of a vector $\boldsymbol{x}$ with a kernel $\boldsymbol{k}$ of size $2n_k + 1$ can be expressed, away from the head and tail of the vector $\boldsymbol{x}$, as
\begin{equation*}
    z_i = (\boldsymbol{x} \star \boldsymbol{k})_i = \sum _{j =- n_k} ^{n_k} x_{i + j} k_{j}.
\end{equation*}
Assuming periodic (circular) boundary conditions and no stride, the matrix $\mathbf{C}$ representing the convolution is the so-called \emph{circulant matrix}. For the sake of clarity, we start with the simple case of one channel input-output layer. Writing the convolution in matrix notation 
$$
\boldsymbol{z} =  \mathbf{C} \boldsymbol{x}
$$
where $\mathbf{C} \in \mathbb{R}^{n_i \times n_i}$ is a circulant matrix defined as
$$
c_{i,j} =
\begin{cases}
k_{j-i}, \quad |j-i|\le n_k \\
0, \qquad \,|j-i|> n_k \\
\end{cases}.
$$
This formulation can be extended to the case where a stride greater than one is applied, while retaining, for simplicity, the assumption of periodic padding. Let $n_{s}$ denote the stride length.
$$
c_{i,j} =
\begin{cases}
k_{j-i}, \quad |j-i|\le n_k \quad  \text{and} \quad  \lfloor (i-1) / n_s  \rfloor = 0 \\
0, \qquad \, \text{elsewhere}\\
\end{cases}.
$$
Given a vector $\boldsymbol{x} \in \mathbb{R}^{n_i }$. Convolving it with $\boldsymbol{k}$ using a stride of $n_s$ produces an output vector $\boldsymbol{z} \in \mathbb{R}^{\tfrac{n_i}{n_s}}$.  

If we replace the convolution operation with a matrix–vector multiplication, the result is a vector of dimension $n_i$, where only $n_i / n_s$ entries are nonzero.
$$
\boldsymbol{z}' =  \mathbf{C} \boldsymbol{x} \in  \mathbb{R}^{n_i }.
$$
Specifically,
$$
z'_i =0 , \quad \text{if} \quad \lfloor (i-1) / n_s  \rfloor \ne 0.
$$
We can build a suitable permutation matrix to push on the bottom the zero rows of the matrix $\mathbf{C}$:
$$
\mathbf{P} \begin{pmatrix}
\mathbf{U} \\
\mathbf{0} \\
\end{pmatrix} \boldsymbol{x}= \boldsymbol{z}',
$$
with $\mathbf{U} \in \mathbb{R}^{n_i/n_s \times n_i}$. We can consequently write:
\begin{equation*}
\label{eq: single LayOrigin}
    \begin{pmatrix}
    \mathbf{U} \\
    \mathbf{0} \\
    \end{pmatrix} \boldsymbol{x}= \mathbf{P}^T \boldsymbol{z}',
\end{equation*}
where $\mathbf{P}^T \boldsymbol{z}'$ has the last $n_i - n_i/n_s$ rows equal to zero. Since $\mathbf{U}$ is a rectangular matrix, it admits the generalized right inverse $\mathbf{U}^+$, which is a matrix so that 
\begin{equation*}
    \mathbf{U} \mathbf{U} ^ + = \mathbb{I} \in \mathbb{R}^{n_{i+1} \times n_{i+1}}.
\end{equation*}
%
Now we introduce the decoder weight matrix defined as $D = \left(\mathbf{U}^+ , \,\, \mathbf{0}^T\right) \mathbf{P}^T$, which gives as output the vector $\boldsymbol{x}'$
$$
\boldsymbol{x}'= D\boldsymbol{z}' = \left(\mathbf{U}^+,\,\, \mathbf{0}\right) \mathbf{P}^T \boldsymbol{z}'.
$$
Notice that in general $\boldsymbol{x}' \ne \boldsymbol{x}$. 
This peculiar choice of the decoder guarantees that the non-zero elements of the vector $\boldsymbol{z}'$ are mapped to themselves through the matrix $E D$, which, expanded, reads as
$$
E D = \mathbf{P} 
\begin{pmatrix}
\mathbf{U} \\
\mathbf{0} \\
\end{pmatrix}  \left(\mathbf{U}^+ , \,\, \mathbf{0}^T\right) \mathbf{P}^T = \mathbf{P}
\begin{pmatrix}
\mathbf{U}\mathbf{U}^+ &  \mathbf{0} \\
\mathbf{0} & \mathbf{0} \\
\end{pmatrix}
\mathbf{P}^T = \mathbf{P}
\begin{pmatrix}
\mathbb{I} &  \mathbf{0} \\
\mathbf{0} & \mathbf{0} \\
\end{pmatrix}
\mathbf{P}^T.
$$
Indeed, clustering in $\bar{\boldsymbol{z}} =( \mathbf{P}^T \boldsymbol{z}')_{1:n_i/n_s} = (\mathbf{P}_{:, 1:n_i/n_s})^T \boldsymbol{z}'$ all the non-zero elements of $\boldsymbol{z}'$, it is clear that $(\mathbf{P}_{:, 1:n_i/n_s})^T E_i D_{i}\mathbf{P}_{:, 1:n_i/n_s}$ maps $\bar{\boldsymbol{z}}$ in $\bar{\boldsymbol{z}}$ and that $(\mathbf{P}_{:, 1:n_i/n_s})^T E_i D_{i}\mathbf{P}_{:, 1:n_i/n_s} = \mathbb{I} \in \mathbb{R}^{n_{i+1} \times n_{i+1}}$.

The extension to convolutional layers with multiple input and output channels is achieved by employing a specific mathematical reformulation that simplifies the formal derivation. To ease the notation, the layer index will be dropped. The size of the $j$-th input channel vector $\boldsymbol{x}_j$ is be denoted by $n$ and $l$ is the number of non zero elements in the $i$-th output channel vector $\boldsymbol{z}_i$. Indeed, let us start by assuming $n_{in}^{ch}, n_{out}^{ch}>1$, and $n_{out}^{ch} \ge n_{in}^{ch}$, where $n_{out}^{ch}$ and $n_{in}^{ch}$ denote the number of output and input channels for a given convolutional layer. The $i$-th output channel vector is the sum of the convolved input channels:
$$
\boldsymbol{z}_i' = \sum _ {j=1} ^{n^{ch}_{in}} \mathbf{C} _{i, j} \boldsymbol{x} _j,
$$
$i = 1,  ..., n_{out}^{ch}$, and $\boldsymbol{z}_i ' \in \mathbb{R} ^ {n} $, $\boldsymbol{x}_j \in \mathbb{R} ^ {n}$.  $\mathbf{C} _{i, j}$ is a single convolution matrix, which represents a convolution operation. Since $j = 1,  ..., n_{in}^{ch}$, a single pass through a multi channels layer involves $n_{in}^{ch} \times n_{out}^{ch}$ convolutions. Following the same notation of the single channel, we can write
$$
\boldsymbol{z}_i' = \sum _ {j=1} ^{n^{ch}_{in}} \mathbf{P} \begin{pmatrix}
\mathbf{U}_{i, j} \\
\mathbf{0} \\
\end{pmatrix}  \boldsymbol{x} _j
\quad \mathrm{and} \; \mathrm{so} \quad 
\mathbf{P} ^ T \boldsymbol{z}_i' = \sum _ {j=1} ^{n^{ch}_{in}} \begin{pmatrix}
\mathbf{U}_{i, j} \\
\mathbf{0} \\
\end{pmatrix}  \boldsymbol{x} _j,
$$
with $\mathbf{U}_{i, j} \in \mathbb{R} ^ {l \times n }$. Assume it is possible to write
$$
\boldsymbol{x}_j' = \sum _ {k=1} ^{n^{ch}_{out}} \left(\mathbf{U}_{j, k}^+ , \,\, \mathbf{0}^T\right) \mathbf{P}^T \boldsymbol{z}_k',
$$
where $\mathbf{U}_{i, j} ^ + \in \mathbb{R} ^ {n \times l}$ then
\begin{align*}
\hat{\boldsymbol{z}}_i &= \sum _ {j=1} ^{n^{ch}_{in}} \mathbf{P} \begin{pmatrix}
\mathbf{U}_{i, j} \\
\mathbf{0} \\
\end{pmatrix}  \boldsymbol{x}' _j  \\
& =  \sum _ {j=1} ^{n^{ch}_{in}} \mathbf{P} \begin{pmatrix}
\mathbf{U}_{i, j} \\
\mathbf{0} \\
\end{pmatrix} \sum _ {k=1} ^{n^{ch}_{out}} \left(\mathbf{U}_{j, k}^+ , \,\, \mathbf{0}^T\right) \mathbf{P}^T \boldsymbol{z}_k' \\
& = \sum _ {j=1} ^{n^{ch}_{in}} \sum _ {k=1} ^{n^{ch}_{out}} \mathbf{P} \begin{pmatrix}
\mathbf{U}_{i, j} \\
\mathbf{0} \\
\end{pmatrix}  \left(\mathbf{U}_{j, k}^+ , \,\, \mathbf{0}^T\right) \mathbf{P}^T \boldsymbol{z}_k' \\
& = \sum _ {j=1} ^{n^{ch}_{in}} \sum _ {k=1} ^{n^{ch}_{out}} \mathbf{P} \begin{pmatrix}
\mathbf{U}_{i, j}\mathbf{U}_{j, k}^+ &  \mathbf{0} \\
\mathbf{0} & \mathbf{0} \\
\end{pmatrix}
\mathbf{P}^T \boldsymbol{z}_k'.
\end{align*}
To have the identity, we have to impose the identity $\hat{\boldsymbol{z}}_i = \boldsymbol{z}_i' $, leading to:
\begin{equation}
\sum _ {j=1} ^{n^{ch}_{in}} \mathbf{U}_{i, j}\mathbf{U}_{j, i}^+ =
\begin{cases}
\mathbb{I} \quad \mathrm{for} \quad k= i\\
\mathbf{0} \quad \mathrm{for} \quad k \ne i
\end{cases}.
\end{equation}
We can get rid of the double sum in all previous relations if we consider the matrices $\mathbf{U}_{j, i}$ and $\mathbf{U}_{j, i} ^ +$ as blocks of a bigger matrix $\overline{\overline{\mathbf{U}}} \in \mathbb{R}^{(l \,n^{ch}_{out})\times (n\, n^{ch}_{in})}$ and
$\overline{\overline{\mathbf{U}}}^+\in \mathbb{R}^{(n\, n^{ch}_{in}) \times (l \, n^{ch}_{out})}$ which are respectively defined as 
$$
\overline{\overline{\mathbf{U}}} =
\begin{pmatrix}
\mathbf{U}_{1,1} & \mathbf{U}_{1,2} & \cdots & \mathbf{U}_{1,n^{ch}_{in}} \\
\mathbf{U}_{2,1} & \mathbf{U}_{2,2} & \cdots & \mathbf{U}_{2,n^{ch}_{in}} \\
\vdots  & \vdots  & \ddots & \vdots  \\
\mathbf{U}_{n^{ch}_{out},1} & \mathbf{U}_{n^{ch}_{out},2} & \cdots & \mathbf{U}_{n^{ch}_{out},n^{ch}_{in}}
\end{pmatrix}
\quad \mathrm{and} \quad 
\overline{\overline{\mathbf{U}}}^+ =
\begin{pmatrix}
\mathbf{U^+}_{1,1} & \mathbf{U^+}_{1,2} & \cdots & \mathbf{U^+}_{1,n^{ch}_{out}} \\
\mathbf{U^+}_{2,1} & \mathbf{U^+}_{2,2} & \cdots & \mathbf{U^+}_{2,n^{ch}_{out}} \\
\vdots  & \vdots  & \ddots & \vdots  \\
\mathbf{U^+}_{n^{ch}_{in},1} & \mathbf{U^+}_{n^{ch}_{in},2} & \cdots & \mathbf{U^+}_{n^{ch}_{in},n^{ch}_{out}}
\end{pmatrix}.
$$
To maintain the consistency it is also necessary to concatenate all the input vectors together, as well as the output vectors. Capital letters denote the vectors obtained by concatenation of the corresponding lowercase letters. For instance:
\begin{equation*}
    \boldsymbol{X} = [\boldsymbol{x}_1^T, \boldsymbol{x}_2^T, ..., \boldsymbol{x}_{n_{in}^{ch}}^T]^ T \in \mathbb{R} ^{n\, n^{ch}_{in}},
\end{equation*}
and
\begin{equation*}
    \boldsymbol{Z}' = [(\mathbf{P}^T\boldsymbol{z'}_1^T)_{1:l}, (\mathbf{P}^T\boldsymbol{z'}_2^T)_{1:l}, ..., (\mathbf{P}^T\boldsymbol{z'}_{n_{out}^{ch}})_{1:l}^T]^ T \in \mathbb{R} ^{l \, n^{ch}_{out}}.
\end{equation*}

Under this formulation, the case of one input and one output channel is a special case of the multiple input-output channels. The condition for the existence of the generalized right inverse ($\overline{\overline{\mathbf{U}}}^+$) ties together the number of input-output channels and the dimension of the input-output vectors. Specifically, $n^{ch}_{out} \,\, l \le n\,\,  n^{ch}_{in}$. 

\subsection{Building the generalised right inverse}
Several strategies can be employed to satisfy the bi-orthogonality condition defined in Equation \ref{eq: bi-ortho condition}. A preliminary approach, following the work of \cite{otto2023learning}, involves parameterizing the generalized inverse as:
\begin{equation}
    \overline{\overline{\mathbf{U^+}}} = \overline{\overline{\mathbf{D}}}(\overline{\overline{\mathbf{U}}} \,\overline{\overline{\mathbf{D}}})^{-1},
\end{equation}
where $\overline{\overline{\mathbf{D}}}$ is a trainable weight matrix. While a judicious construction of $\overline{\overline{\mathbf{D}}}$ allows for a low-rank parameterization of the decoder, numerical experiments indicated that this formulation suffers from poor training convergence. Specifically, the loss function frequently stagnated, or the model reached a pathological state where the weights failed to represent a physically meaningful manifold.

To overcome these limitations, we propose an alternative formulation based on orthogonal projectors. In this framework, the decoding stage is not restricted to a linear mapping but is instead decomposed into a linear component and a non-linear term:
\begin{equation}
\label{eq: projecion layer}
    \boldsymbol{Z}'_o = D(\boldsymbol{Z}'_i) =  \underbrace{\overline{\overline{\mathbf{U}}} ^{P} \boldsymbol{Z}'_i }_{\text{Linear}}+  \underbrace{(I - \overline{\overline{\mathbf{U}}} ^{P} \overline{\overline{\mathbf{U}}})\boldsymbol{f}(\boldsymbol{Z}'_i)}_{\text{Non Linear}}
\end{equation}
where $\overline{\overline{\mathbf{U}}} ^{P}$ denotes the Moore-Penrose pseudoinverse. Here, the non-linear contribution is governed by a function $\boldsymbol{f}$, which in our implementation consists of classical deconvolutional layers. This architecture offers several theoretical advantages:

\begin{itemize}
    \item If $\boldsymbol{f}$ is the null function, the decoder becomes parameter-free, reducing to a Symmetric Orthogonal Autoencoder (SOAE) as discussed in \cite{brivio2025deep}.
    \item When $\boldsymbol{f}(\boldsymbol{z}'_i) \neq \boldsymbol{0}$, it provides the decoder with additional degrees of freedom to capture complex features of the physical solution without violating the bi-orthogonality constraint.
    \item The non-linear term is designed to lie within the null space of the encoder. Consequently, the identity $E \circ D = id$ is preserved regardless of the complexity of $\boldsymbol{f}$. Keeping in mind that $\overline{\overline{\mathbf{U}}} \overline{\overline{\mathbf{U}}} ^{P} = \mathbb{I}$, this is formally demonstrated by:
    \begin{gather*}
        \boldsymbol{Z} = E(\boldsymbol{X}) = \overline{\overline{\mathbf{U}}} \boldsymbol{X} \\
        \boldsymbol{X}' = D(\boldsymbol{Z} )
    \end{gather*}
    \begin{equation*}
                \boldsymbol{Z}'  = E(\boldsymbol{X}') = E( D(\boldsymbol{Z} )) = 
          \underbrace{\overline{\overline{\mathbf{U}}} \overline{\overline{\mathbf{U}}} ^{P} }_{\mathbf{I}} \boldsymbol{Z} +  \underbrace{\overline{\overline{\mathbf{U}}} (I - \overline{\overline{\mathbf{U}}} ^{P} \overline{\overline{\mathbf{U}}})}_{\boldsymbol{0}} \boldsymbol{f}(\boldsymbol{Z}) = \boldsymbol{Z}.
    \end{equation*}
    This ensures that any point decoded from the latent space is perfectly preserved upon re-encoding, providing a consistent and stable coordinate system for the reduced-order model.
\end{itemize}

\subsection{Activation functions}
As previously discussed, the selection of appropriate activation functions is critical for maintaining a globally representation consistent AE. To ensure that the composition of the encoder $E$ and decoder $D$ yields the identity mapping across multiple layers, each activation function must be invertible. Let $\sigma$ denote a Bi-Lipschitz, invertible activation function that operates element-wise. The global consistency of the network is preserved through the following chain of compositions:
\begin{equation}
    E \circ D =  E_1 \circ \sigma \circ \cdots \circ \sigma \circ E_n  \circ D_n \circ \sigma^{-1} \circ \cdots \circ \sigma^{-1} D_1 =  E_1 \circ D_1 = id.
\end{equation}
It is clear from the expanded chain that the overall composition results into the identity. The invertible activation functions already existing in the literature are the LeakyReLU \cite{brivio2025deep} and the Hyperbolic activation \cite{brivio2025deep, otto2023learning}. We enlarge this collection by proposing a new invertible and Bi-Lipschitz activation function resembling the hyperbolic tangent that we name \emph{Rational} activation function:
\begin{equation}
    \label{eq: rational AF}
    y = \sigma_\alpha(x) = \alpha \left( x + \frac{x}{1+x^2}\right),
\end{equation}
where $\alpha$ accounts as hyperparameter. The inverse is analytic and can be written as:
\begin{equation*}
    x = \sigma ^{-1}_\alpha(y) = \left( - \frac{q}{2} + \sqrt{\Delta} \right)^{1/3} + \left( - \frac{q}{2} - \sqrt{\Delta} \right)^{1/3},
\end{equation*}
with
\begin{align*}
t&=\frac{y}{\alpha},           &  p &=2 -\frac{t^2}{3},             \\
q&=-\frac{2t^3+9t}{27},        &  \Delta &= \left(\frac{q}{2}\right)^2 +  \left(\frac{p}{3}\right)^2 .\\    
\end{align*}

\begin{figure}[h!]
    \centering
    \includegraphics[width=0.75\linewidth]{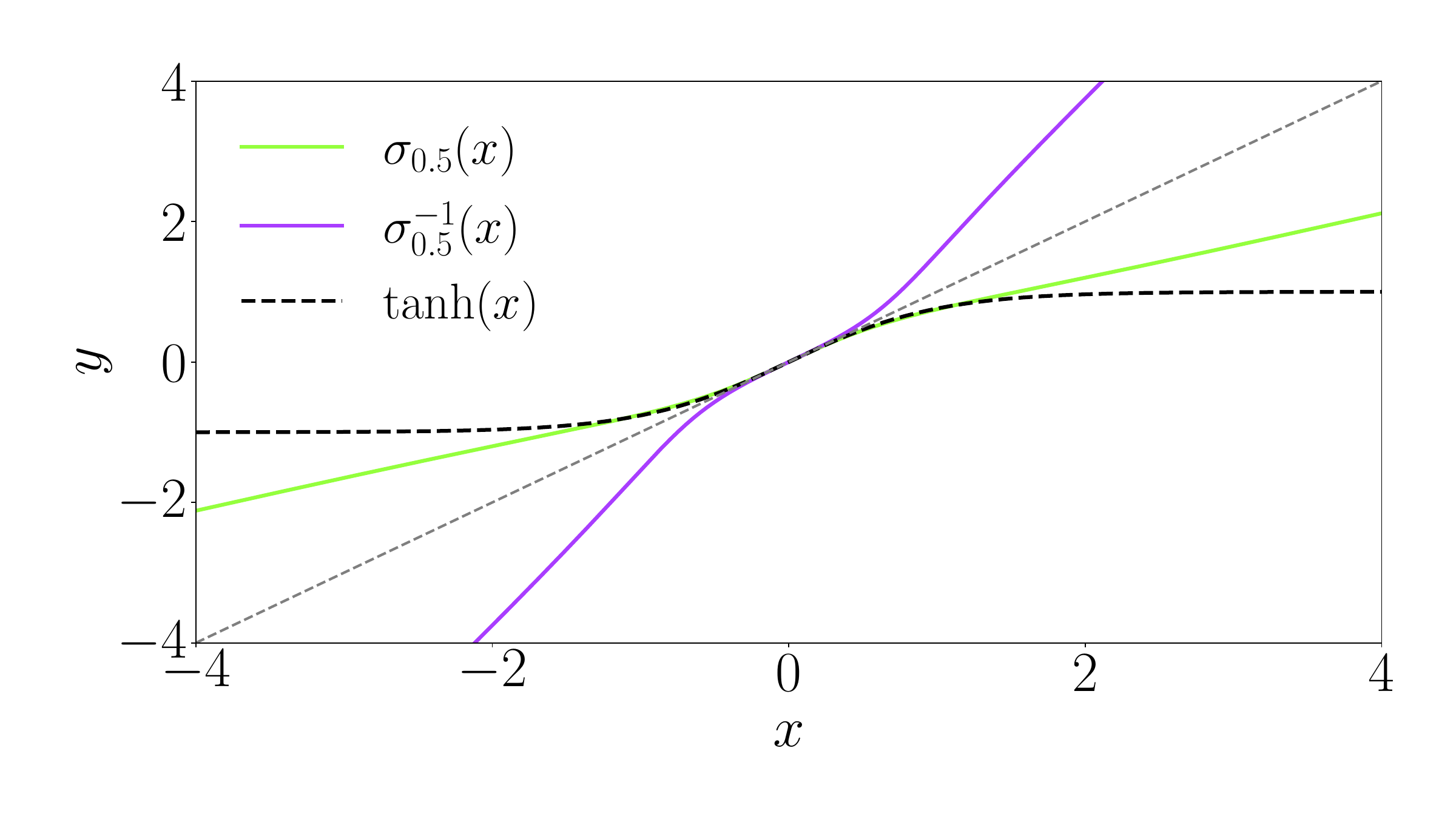}
    \caption{Rational activation function $\sigma_\alpha$, its inverse $\sigma_\alpha ^{-1}$, and Tanh activation function. The dashed gray line denotes the bisector of the first and third quadrants.}
    \label{fig:AF}
\end{figure}

The upper and lower Lipschitz constants of the rational activation function are
\begin{equation*}
    L = 2\alpha \quad \mathrm{and} \quad \eta = \alpha \frac{7}{8}
\end{equation*}
such that 
\begin{equation*}
    \eta\|x - y \| \le \|\sigma_\alpha(x) - \sigma_\alpha(y)\| \le L  \| x - y \|, \qquad \forall x, y \in \mathbb{R}.
\end{equation*}
We conduct all the numerical experiments with $\alpha = 0.5$. A plot of the rational activation function and its inverse, for the given value of the hyperparameter $\alpha$ is given in Figure \ref{fig:AF}.

\section{Data-driven approach}
\label{sec: Data-Driven}
Having introduced the enhanced version of convolutional autoencoders, which are distinguished by their advantageous parameterization properties, we can now employ them for manifold learning tasks to obtain a reduced representation of high-dimensional data. Additionally, we can leverage the universal approximation capabilities of feedforward neural networks to model the function that defines the latent flow described in Equation \ref{eq: latent equation}.

To construct the reduced representation of the dynamical systems governed by Equation \ref{eq: FOM equation}, we integrate the aforementioned architectures within a decoupled two-stage training framework. In the first stage, the autoencoder is trained on the available snapshot data to minimize reconstruction error. Once the AE has converged, its weights are frozen, and the encoder is utilized to project the high-fidelity trajectories into the latent space. In the second stage, these latent trajectories serve as the ground truth for training the LF-net, ensuring that the learned flow accurately captures the time evolution of the reduced system.

\subsection{AutoEncoders}
\label{section: AE}
The maps $\mathbf{x}_\alpha$ are not known a priori, but an approximation of them can be learned by means of a data-driven approach. Indeed, let us assume to have access to a certain number of realizations of the solution of Equation \ref{eq: FOM equation} for a given set of parameters $\boldsymbol{P} \subset \mathcal{P}$. We call this set of snapshots
\begin{equation*}
    \mathcal{D} =\{ \mathbf{u}_h(t; \boldsymbol{\mu}) | t \in \Theta , \, \mu \in \boldsymbol{P} \},
\end{equation*}
where $\Theta$ is a set of discrete time instants. By exploiting the manifold hypothesis, we assume that all the points in the dataset live in a low dimensional manifold, which we may want to learn by using a suitable neural network architecture such as an AE. Specifically, the decoder $D: \mathbb{R}^l \rightarrow \mathbb{R}^{N_h}$ tends to approximate the map $\mathbf{x}_\alpha$ introduced in Equation \ref{eq: exact decoder}, while the encoder $E: \mathbb{R}^{N_h} \rightarrow \mathbb{R}^{l}$ is an approximation of its inverse. The AutoEncoder is defined as the composition:
\begin{equation}
    AE = D \circ E : \mathbb{R}^{N_h} \rightarrow \mathbb{R}^{N_h}.
\end{equation}
Training the AE on the dataset $\mathcal{D}$ consists in finding the maps $E$ and $D$ which minimize the reconstruction error. Namely, 
\begin{equation}
\label{eq: minimization formulation}
    \min_{E, D} \text{MSE} (\boldsymbol{u}, \,\, AE(\boldsymbol{u})), \,\, \boldsymbol{u} \in \mathcal{D}.
\end{equation}
The combined action of decoder and encoder provide an approximation of the solution manifold:
\begin{equation*}
    \widetilde{S}_h ^n\approx \{AE(u) | \, u \in S_h\}.
\end{equation*}
Since $\mathbf{u}_h \simeq \widetilde{\mathbf{u}} \in \widetilde{S}_h $ it is possible to assume that $\mathbf{u}_h $ satisfies Equation \ref{eq: FOM equation} as well, leading to:
\begin{equation}
\label{eq: approx FOM eq.}
    \left\{\begin{array}{l}
\dot{\widetilde{\mathbf{u}}}(\boldsymbol{z}(t ; \boldsymbol{\mu}))=\mathbf{f}_h\left(t, \widetilde{\mathbf{u}}(\boldsymbol{z}(t ; \boldsymbol{\mu})) ; \boldsymbol{\mu}\right), \quad t \in\left(t_0, T\right], \\
\widetilde{\mathbf{u}}\left(\boldsymbol{z}(t_0 ; \boldsymbol{\mu})\right)=\widetilde{\mathbf{u}}_{0}(\boldsymbol{\mu}).
\end{array}\right.    
\end{equation}
After some careful manipulations, the latent governing equations produced by the AEs maps read: 
\begin{equation}
    \label{eq: latent equation AE}
    \dot{\boldsymbol{z}}(\boldsymbol{z}, t ; \boldsymbol{\mu})) = (dEdD)^{-1} dE \mathbf{f}_h\left(t, \widetilde{\mathbf{u}}(\boldsymbol{z}(t ; \boldsymbol{\mu})) ; \boldsymbol{\mu}\right).
\end{equation}
Before proceeding any further, several considerations are warranted. First, a striking similarity can be observed between Equation \ref{eq: latent equation AE}, derived from true manifold parametrizations, and Equation \ref{eq: latent equation}. Although these two equations differ only by the additional term $(dE \, dD)^{-1}$ present in Equation \ref{eq: latent equation AE}, this term vanishes if the AE belongs to the family of symmetric AEs, thereby reducing Equation \ref{eq: latent equation AE} to Equation \ref{eq: latent equation}. Despite the fact that Equation \ref{eq: latent equation AE} is globally low-dimensional and allows for the autonomous evolution of the dynamics, the temporal integration of such a system would offer no computational advantage, as evaluating the RHS of Equation \ref{eq: latent equation} requires intermediate computation of high-dimensional vectors. Any type of nonlinearity appearing in the flow $\mathbf{f}_h\!\left(t, \widetilde{\mathbf{u}}(\boldsymbol{z}(t;\boldsymbol{\mu}))\right)$ or in the model order reduction technique itself, in fact, leads to such a caveat. To circumvent this issue, an appropriate neural network architecture can be introduced to approximate the latent flow, thereby providing the efficiency required for a high-performance ROM.  It is worth stressing that the proposed approach is therefore particularly well suited for non-linear parametric PDEs. First, non-linear PDEs often suffer from a slow decay of the Kolmogorov $n$-width (\eg, turbulence, wave propagation), thus highlighting the need for non-linear model order reduction techniques such as AEs. Secondly, the need for a dedicated neural network to advance the latent dynamics in time is intrinsically related to the non-linearity of the FOM equations themselves.

\subsection{Latent Flow Network (LF-net)}
\label{subsec:LF_net}
\begin{figure}[h!]
    \centering
    \includegraphics[width=0.3\linewidth]{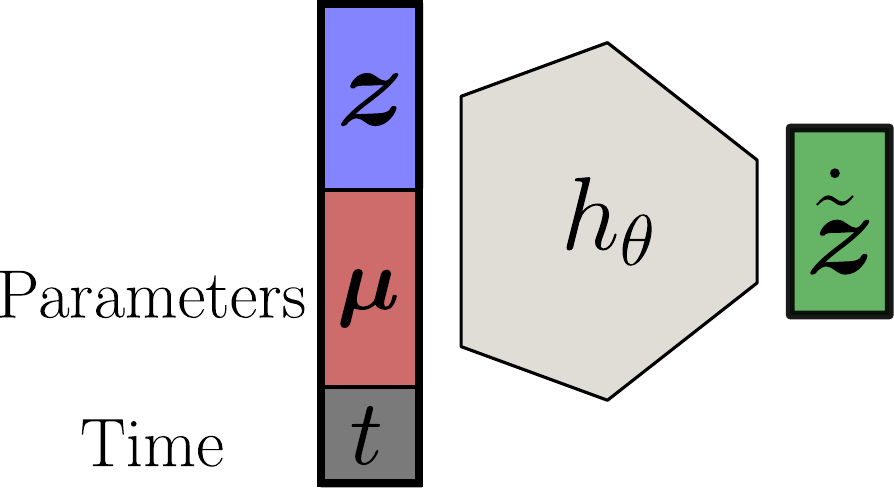}
    \caption{Sketch of the LF net.}
    \label{fig:LF net}
\end{figure}

The latent flow is approximated by a Feed-Forward Neural Network, denoted as $h_\theta$. According to the dynamics established in Equation \ref{eq: latent equation}, the flow typically depends on the current latent state $\boldsymbol{z}$, the parameter vector $\boldsymbol{\mu}$, and, in non-autonomous cases, the explicit time $t$. To account for these dependencies, these variables are concatenated into a single vector, which is given as input to the network, as illustrated in Figure \ref{fig:LF net}.
Specifically, the latent state is obtained through the encoding mapping $\boldsymbol{z} = E(\mathbf{u}_h(t; \boldsymbol{\mu}))$. The ground-truth latent velocity $\dot{\boldsymbol{z}}$ is derived by computing the differential of the encoder with respect to the high-fidelity state $\mathbf{u}_h$, evaluated at $\mathbf{u}_h(t; \boldsymbol{\mu})$, and applying to the high-fidelity dynamics $\mathbf{f}_h$ as follows:
\begin{equation}
    \dot{\boldsymbol{z}}(\boldsymbol{z}, t; \boldsymbol{\mu}) =dE\left(\mathbf{u}_h(t; \boldsymbol{\mu})\right) \cdot \mathbf{f}_h\left(t, \mathbf{u}_h(t ; \boldsymbol{\mu}) ; \boldsymbol{\mu}\right).
\end{equation}
The LF-net is trained by minimizing a Mean Squared Error (MSE) loss function, which penalizes the discrepancy between the network's prediction and the projected high-fidelity dynamics:
\begin{equation}
    \min_{h_\theta} \text{MSE} \{h_\theta(E(\boldsymbol{u}_h), \boldsymbol{\mu}, t) , dE \mathbf{f}_h(t, \boldsymbol{u}_h;  \boldsymbol{\mu}) \} \,\, \boldsymbol{u}_h \in \mathcal{D}.
\end{equation}
\subsection{ROM Deployment}
\label{subsec:ROM_deployment}

Following the offline training of both the AE and the LF-net, the ROM can be deployed for efficient online predictions. The deployment workflow consists of two main steps: initial condition projection and latent trajectory evolution. First, the high-fidelity initial condition $\mathbf{u}_{0, h}(\boldsymbol{\mu})$ is projected into the latent space using the learned encoder $E$:
\begin{equation}
    \boldsymbol{z}(0 ; \boldsymbol{\mu}) = E(\mathbf{u}_{0, h}(\boldsymbol{\mu})).
\end{equation}
Subsequently, the latent state is evolved in time by integrating the learned latent flow $h_\theta$. This is achieved using a suitable numerical integration scheme (e.g., explicit or implicit Euler, Heun, or higher-order Runge-Kutta methods):
\begin{equation}
    \label{eq: ROM solution}
    \widetilde{\boldsymbol{z}}(t; \boldsymbol{\mu}) = \boldsymbol{z}(0 ; \boldsymbol{\mu}) + \int_{0}^{t} h_\theta(\widetilde{\boldsymbol{z}}(\tau), \boldsymbol{\mu}, \tau) \, d\tau.
\end{equation}
%
In general, integrating \( h_\theta \) does not equate to solve Equation \ref{eq: latent equation AE}. The primary reason for this is that \( h_\theta \) is designed to approximate the projection of the tangent vector within the latent tangent space learned by the encoder. However, in the case of a general autoencoder, this projection does not adequately represent the latent space governing equation, as it neglects the term \( (dEdD)^{-1} \), accounting for the representation-inconsistent parameterization that has been learned. Conversely, if the encoder and decoder are chosen from families of functions such that \( E \circ D = id \) (i.e., a representation-consistent formulation), then Equation \ref{eq: ROM solution} is the solution of Equation \ref{eq: latent equation AE}.

\section{Numerical experiments}
\label{sec: Numerical experiments}

\begin{figure}[h!]
\centering
\begin{subfigure}{0.3\linewidth}
    \centering
    \includegraphics[width=0.9\linewidth]{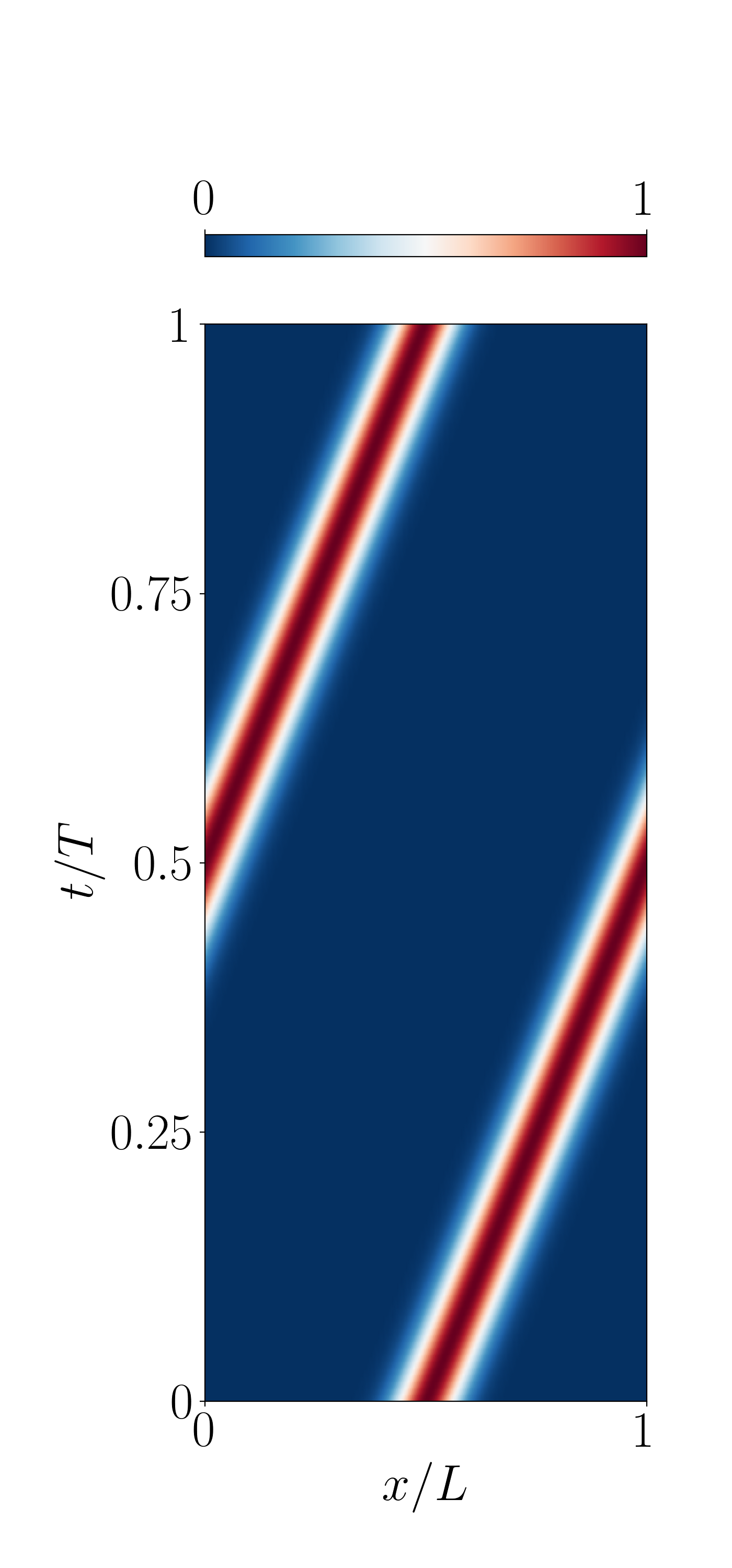}

\end{subfigure}
\begin{subfigure}{0.3\linewidth}
    \centering
    \includegraphics[width=0.9\linewidth]{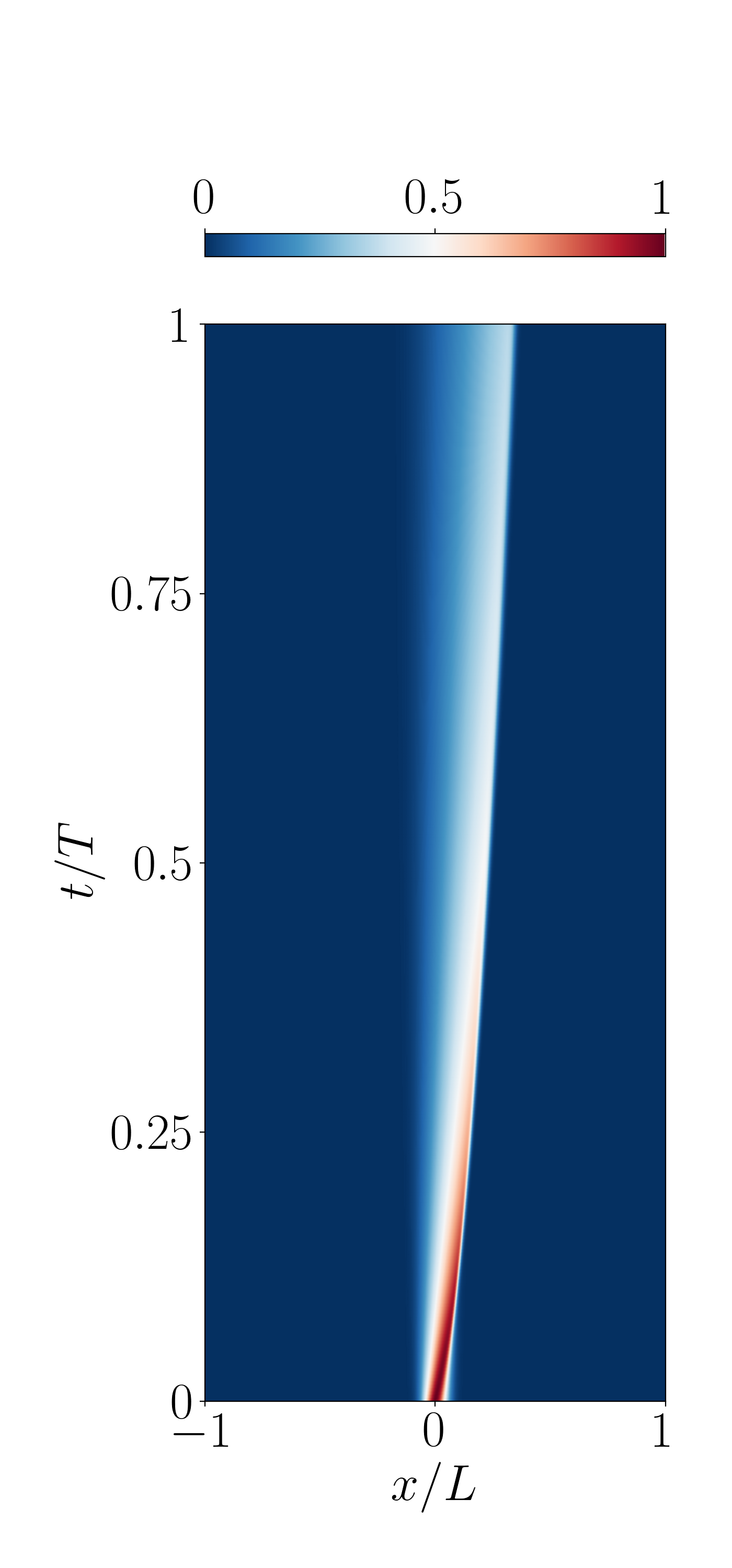}

\end{subfigure}
\begin{subfigure}{0.3\linewidth}
    \centering
    \includegraphics[width=0.9\linewidth]{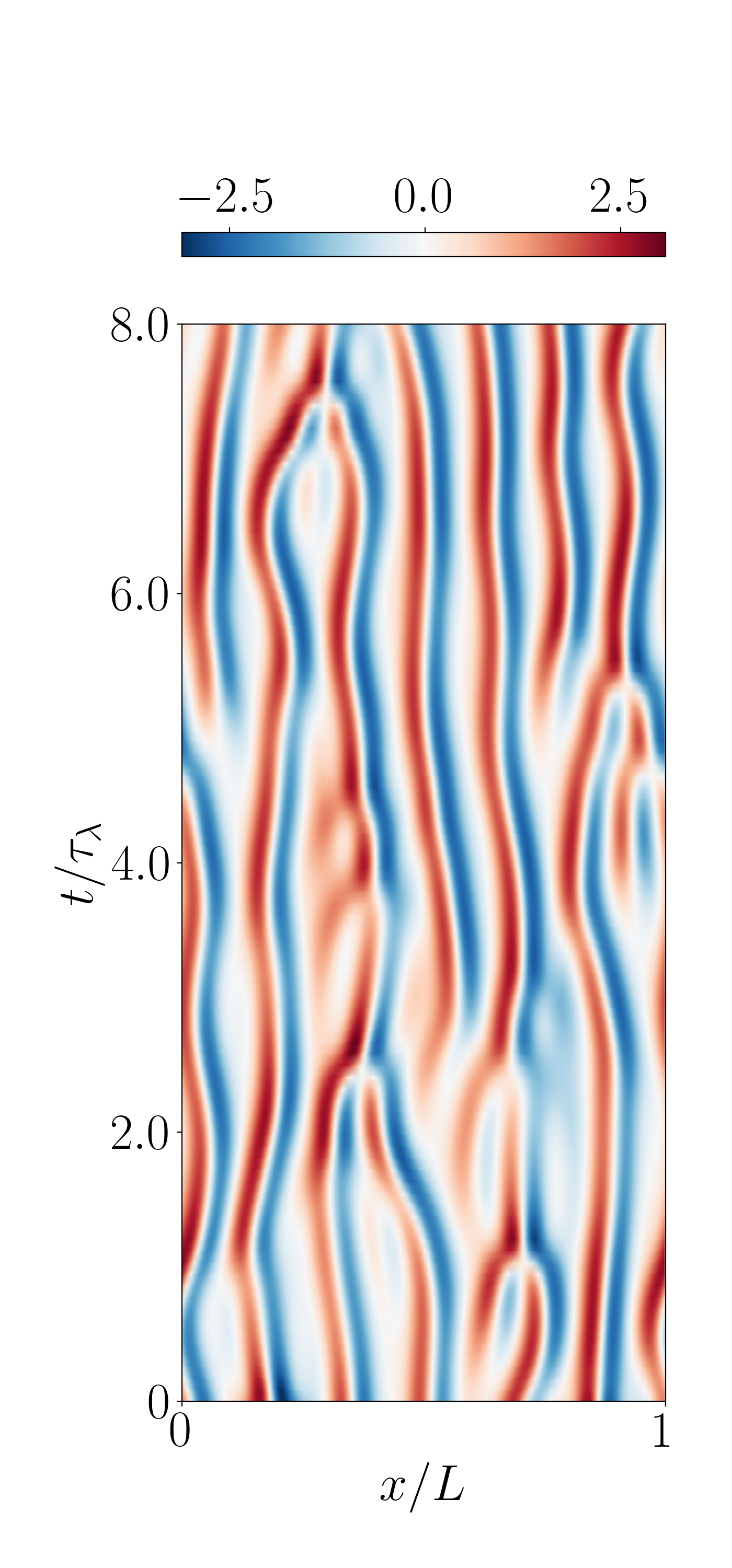}

\end{subfigure}

\caption{Solution of Linear Advection equation (left), Viscous Burgers (middle) and Kuramoto-Sivanshisky equation (right).}
\label{fig: Physical Systems}
\end{figure}

To evaluate the proposed Convolutional Symmetric Autoencoders within a Model Order Reduction framework, the pipeline detailed in Section \ref{sec: Data-Driven} is applied to three canonical one-dimensional dynamical systems, whose spatio-temporal solutions are depicted in Figure \ref{fig: Physical Systems}. Each system was selected to represent a distinct challenge, thereby demonstrating the versatility of the architecture across different levels of complexity. The assessment focuses on both reconstruction accuracy, computational cost and latent space stability. 
We define latent stability as the capacity of the reduced model to produce solutions that are physically consistent. In deterministic, non-chaotic systems, this means that the latent trajectories remain close to the ground truth. However, in chaotic systems, this condition becomes impractical because, by nature, trajectories that initially start close will eventually diverge. In these cases, stability is characterized by the ability of the latent trajectories to remain within the invariant set, which is represented by the projection of the high-dimensional attractor into the latent space.
The analysis begins with Linear Advection, a benchmark for advection-dominated problems where traditional linear reduction techniques, such as POD, typically underperform. This is followed by the parametric Viscous Burgers’ equation, where the viscosity serves as the primary parameter. Finally, the model is tested against the Kuramoto-Sivashinsky equation, representing a chaotic system.

\subsection{Linear Advection}
\label{sec: linear_advection}
The one-dimensional  Linear Advection equation describes the transport of a passive scalar field within a spatial domain of size $L$ over a temporal horizon $T$. The governing equations, assuming periodic boundary conditions, are defined as:
\begin{equation}
\label{eq: Advection equation}
\begin{cases}
\frac{\partial  u}{\partial t} + v \frac{\partial u} {\partial x} = 0, \quad x \in [0, L], t \in [0, T]\\
u(x, 0) = u_0(x) \\
u(0, t) = u(L, t)  
\end{cases}
\end{equation}
where $v$ denotes the constant transport velocity. Despite its mathematical simplicity and linearity, this equation remains a significant challenge for classical Model Order Reduction techniques. In advection-dominated regimes, the Kolmogorov $n$-width exhibits a slow decay rate, meaning that linear subspace methods, such as POD, require an excessively large number of modes to achieve accurate compression. Recently it was shown that the decay of the Kolmogorov $n$-width is also largely affected by the regularity and smoothness of initial and boundary conditions \cite{arbes2025kolmogorov}.

Equation \ref{eq: Advection equation} is solved numerically using a second-order central finite difference scheme as spatial discretization and a $4^{\mathrm{th}}$-order five stages Runge-Kutta method for the time integration. For this benchmark, the physical parameters, specifically the advection velocity $v$, the simulation time $T$, and the domain length $L$, were set to unity.
The initial condition $u_0(x)$ was defined as a Gaussian centered in the middle of the domain with a normalized standard deviation of $\sigma/L = 0.05$. The spatial interval $[0, L]$ was partitioned into $N_h + 1$ equidistant points with a grid spacing of $h = L/N_h$. Let $\boldsymbol{u}_h \in \mathbb{R}^{N_h}$ represent the discrete solution vector, where $u_i \approx u(x_i, t)$ for $x_i = ih$ and $i = [0, 1, ..., N_h-1]$. The resulting semi-discrete system is given by:
\begin{equation}
\label{eq: semi_discrete_advection}
\frac{d \boldsymbol{u}}{d t} + \mathbf{A} \boldsymbol{u} = 0.
\end{equation}
The dataset used for training and testing consists of snapshots derived from these numerical solutions. To ensure the generalization of the model, $80\%$ of the data was allocated for training, while the remaining $20\%$ was reserved for testing and validation.

\subsection{Viscous Burgers Equation}

The Viscous Burgers' equation is a fundamental one-dimensional non-linear parametric PDE that serves as a benchmark for modeling advection-diffusion phenomena. In this study, we consider the equation in its conservative form under periodic boundary conditions:
\begin{equation*}
\begin{cases}
\frac{\partial  u}{\partial t} +  \frac{\partial  } {\partial x} \left( \frac{u ^ 2}{2}  \right)  = \nu \frac{\partial ^2 u} {\partial x ^2}  \quad x \in [-L, L],\, t \in [0, T],\, \nu \in (0.001, 0.1],\\
u(x, 0) = u_0(x) \\
\text{Periodic B.C.s}
\end{cases}
\end{equation*}
where the viscosity $\nu \in  (0.001, 0.1]$ acts as the system parameter. The non-linear convective flux is primarily responsible for the development of steep gradients and shock-like discontinuities. Such localized features are notoriously difficult for global linear bases to capture. This limitation often manifests itself as spurious oscillations around the discontinuity \cite{fresca2021comprehensive, maulik2021reduced}.

The governing equation is discretized using a FD framework with $L=1$ and $T=1$. Spatial derivatives are approximated using a second-order central difference scheme on a uniform grid with $N_h + 1$ discretization points, resulting in a uniform mesh size $h = 2L/N_h$. For temporal integration, an explicit Euler scheme is employed. 

To construct the ROM, the parameter interval for $\nu$ is sampled at $15$ logarithmically spaced points. This sampling strategy ensures that the model is exposed to a wide range of physical regimes, from diffusion-dominated to advection-dominated flows. Out of the $15$ generated trajectories, $12$ are used for the training phase, while the remaining 3 are reserved for testing the model.

\subsection{Kuramoto-Sivashinsky Equation}
\label{sec:ks_equation}

The final test case considers the Kuramoto-Sivashinsky (KS) equation, which introduces a significantly higher level of complexity due to its chaotic dynamics. Originally developed to model the instabilities of laminar flame fronts, the KS equation is one of the simplest PDEs known to exhibit high-dimensional chaos \cite{huntsman2023towards}. The governing equation, under periodic boundary conditions, is given by:
\begin{equation}
\label{eq:ks_equation}
\begin{cases}
\frac{\partial u}{\partial t} + u \frac{\partial u}{\partial x} + \frac{\partial^2 u}{\partial x^2} + \nu \frac{\partial^4 u}{\partial x^4} = 0, & x \in [0, L], \ t \in [0, T] \\
u(x, 0) = u_0(x) \\
u(0, t) = u(L, t) .
\end{cases}
\end{equation}
The qualitative behavior of the solution, ranging from a bursting regime to fully developed chaos, is strictly dependent on the viscosity $\nu$ and the domain length $L$ \cite{COLANERA2025118393}. In this study, we focus on the chaotic regime, which is readily achieved by selecting a sufficiently large domain $L$ \cite{edson2019lyapunov}. Indeed the number of modes that growth exponentially increases linearly with $L$. In the chaotic regime, the nonlinear term plays a stabilizing role by facilitating an energy cascade, effectively redistributing energy from unstable low-frequency modes to damped high-frequency modes \cite{baez2022kuramotosivashinskyequation}.

The system is solved numerically using a pseudo-spectral Fourier method for spatial discretization, coupled with an Exponential Time Differencing Fourth-order Runge-Kutta (ETDRK4) scheme for time marching, as detailed in \cite{kassam2005fourth}. 
A more comprehensive discussion on the properties of chaotic systems and their reduced-order modeling is provided in Appendix \ref{sec: chaotic systems}.

\section{Results}
\label{sec: Results}
\subsection{LA equation}

The study of the Linear Advection problem is motivated by two primary factors. First, for non-linear reduction techniques, it facilitates the generation of a two-dimensional latent space, which lends itself to straightforward visualization and interpretability. Second, it provides a rigorous benchmark for evaluating the relative performance of three distinct methodologies: POD, classical Convolutional Autoencoders, and their symmetric counterparts.

To construct the POD-based ROM, the numerical solutions are organized into a snapshot matrix $\boldsymbol{X} = [\boldsymbol{u}_1, \boldsymbol{u}_2, \dots, \boldsymbol{u}_{N_T}] \in \mathbb{R}^{N_h \times N_T}$, where $\boldsymbol{u}_i$ denotes the discrete solution to Equation \ref{eq: semi_discrete_advection} at time $t_i$. By computing the Singular Value Decomposition (SVD) of $\boldsymbol{X}$, the POD modes $\boldsymbol{\varphi}_i$ (the left singular vectors) are extracted to span a reduced-order subspace. The governing equation is then projected onto this subspace to derive the evolution equations for the POD coefficients, defined as $c_i = \langle \boldsymbol{\varphi}_i, \boldsymbol{u} \rangle$:
\begin{equation}
\label{eq:pod_evolution}
\frac{d \mathbf{c}}{d t} + \boldsymbol{\Phi}^T \mathbf{A} \boldsymbol{\Phi} \mathbf{c} = 0
\end{equation}
where $\boldsymbol{\Phi} = [\boldsymbol{\varphi}_1, \boldsymbol{\varphi}_2, \dots, \boldsymbol{\varphi}_l]$ represents the basis truncated at dimension $l$. According to the Eckart–Young theorem, the truncated SVD provides the optimal $l$-rank approximation of $\boldsymbol{X}$ in the Frobenius norm. Furthermore, the total reconstruction error is bounded by the sum of the squares of the discarded singular values $\sigma_i$ from $i = l+1$ to $r$, where $r = \min\{N_h, N_T\}$ is the rank of the snapshot matrix \cite{brunton2022data}.

\begin{figure}[h!]
    \centering
    \includegraphics[width=0.75\linewidth]{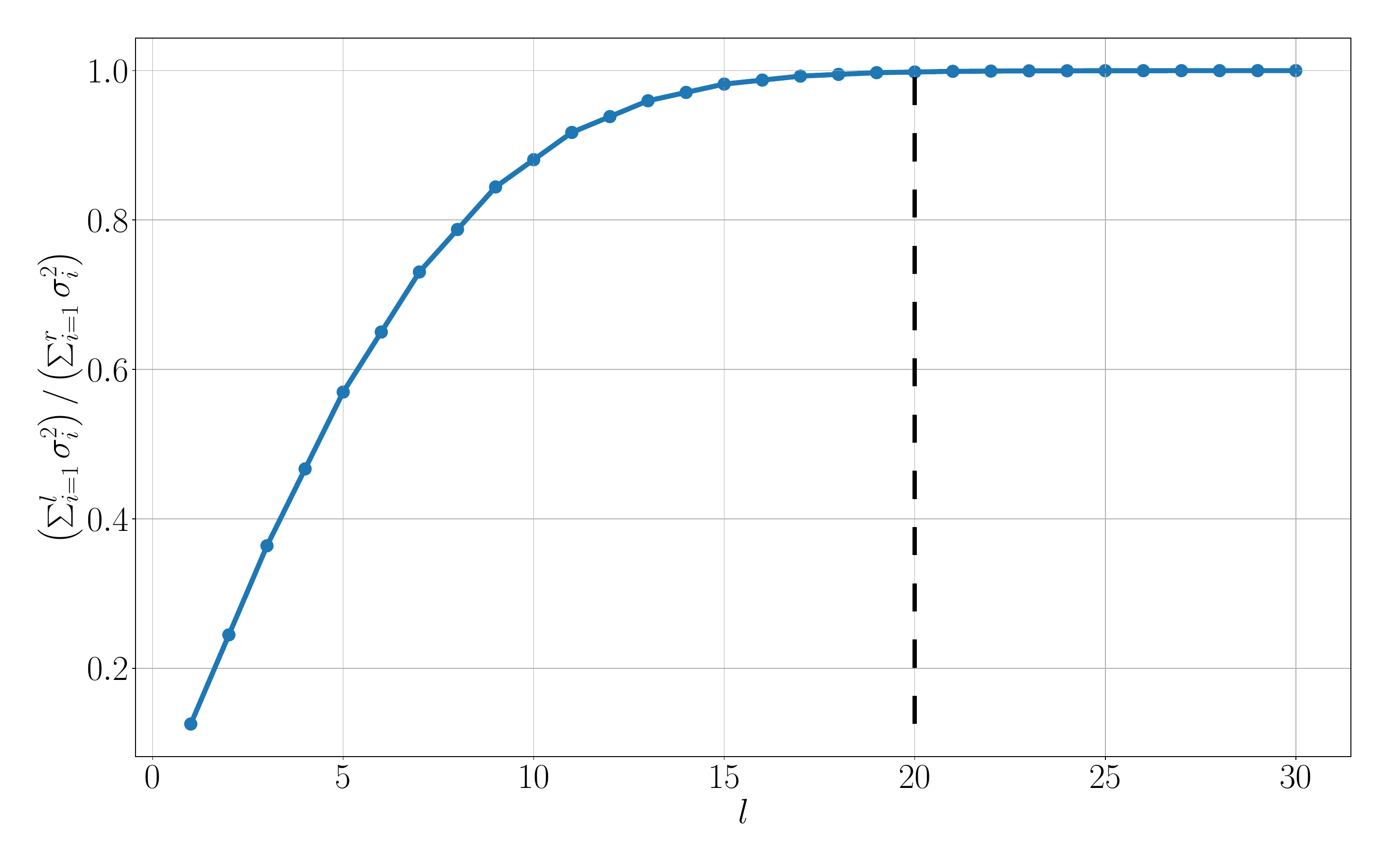}
    \caption{Cumulative sum of the squared singular values of the snapshot matrix $\boldsymbol{X} = [\boldsymbol{u}_1, \boldsymbol{u}_2, \dots, \boldsymbol{u}_{N_T}]\in \mathbb{R}^{N_h \times N_t }$}
    \label{fig: Cumulative Sum}
\end{figure}

To estimate the truncated dimension $l$ required for an accurate POD-ROM in the unit domain problem described in Section \ref{sec: linear_advection}, we analyze the cumulative energy carried by the singular values of $\boldsymbol{X}$. As illustrated in Figure \ref{fig: Cumulative Sum}, the cumulative sum of the squared singular values suggests that approximately $l \simeq 20$ modes are sufficient to capture the dominant ($99.82\%$) information within the dataset.
To quantify the performance of the ROMs, we define the relative $L_2$ reconstruction error as:
\begin{equation}
\label{eq: epsilon u}
\epsilon_{\boldsymbol{u}}(t) = \frac{\| \boldsymbol{u}(t) - \tilde{\boldsymbol{u}}(t) \|_{L_2}}{\| \boldsymbol{u}(0) \|_{L_2}}
\end{equation}
where $\tilde{\boldsymbol{u}}(t)$ is the solution predicted by the ROM. For the POD-ROM, the reconstruction is given by $\tilde{\boldsymbol{u}}(t) = \boldsymbol{\Phi} \mathbf{c}(t)$. As illustrated in Figure \ref{fig:  Linear Advection POD error}, a POD basis of $l=20$ is required to maintain the relative error below $1\%$ throughout the simulation. Reducing the number of modes to $l < 10$ results in significantly higher errors, underscoring the limitations of linear reduction techniques even for relatively simple advection problems.
\begin{figure}[h!]
    \centering
    \includegraphics[width=0.75\linewidth]{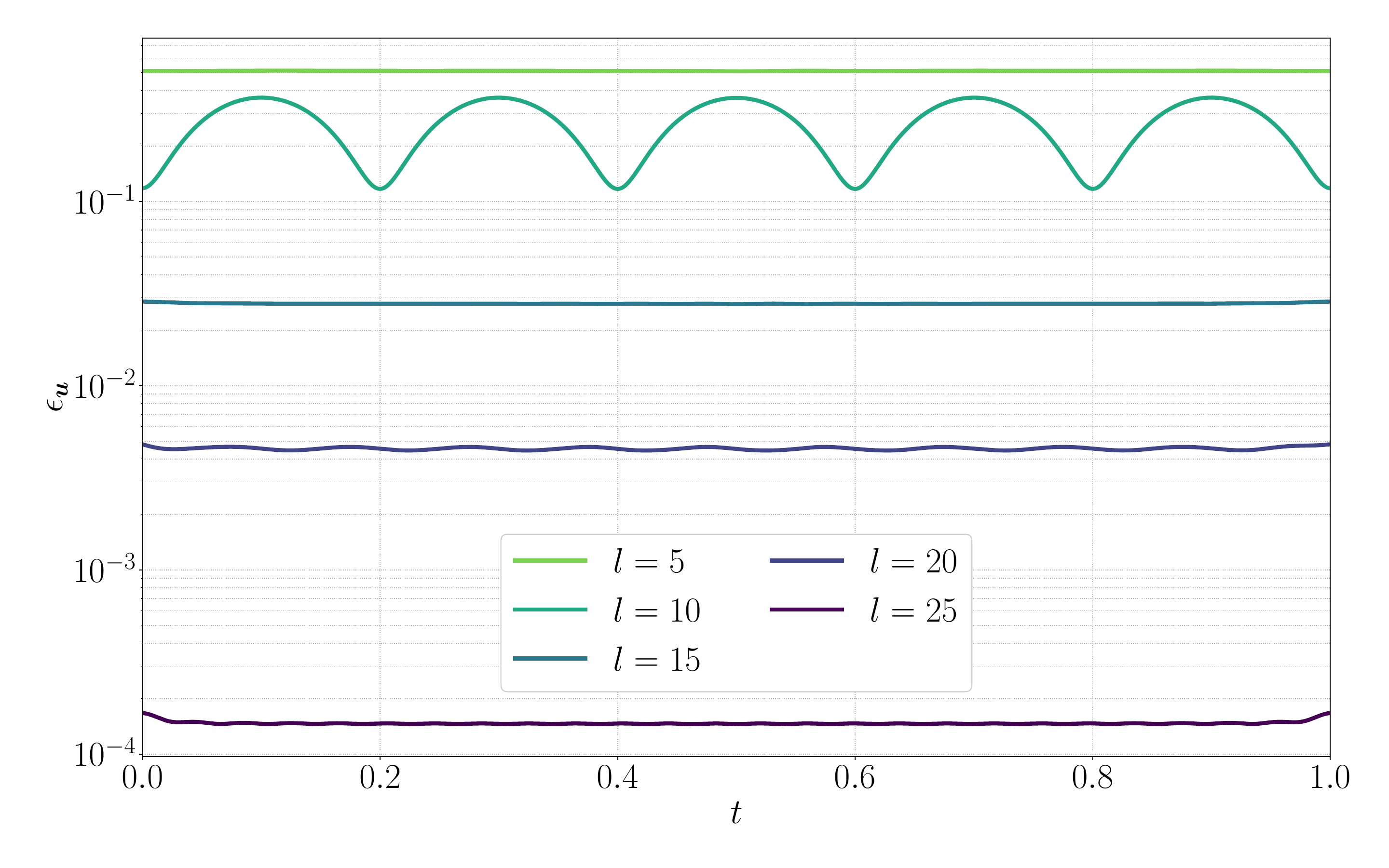}
    \caption{Error $\epsilon _ {\boldsymbol{u}}$ defined in Equation \ref{eq: epsilon u} with an increasing number of modes $l$, as function of time $t$.}
    \label{fig: Linear Advection POD error}
\end{figure}
The limitations of POD motivate the use of Autoencoders for non-linear manifold learning. Given that the solution to Equation \ref{eq: Advection equation} is periodic and governed by a single parameter (velocity), the FOM solution traces a closed curve embedded in $\mathbb{R}^{N_h}$. While this manifold is locally one-dimensional, its periodic nature prevents a global, injective parametrization into $\mathbb{R}$ \cite{floryan2022data}. Consequently, a minimum latent dimension of $l=2$ is required for the AE to capture the topological feature of the solution.
\begin{figure}[h!]
     \includegraphics[width=0.7\linewidth]{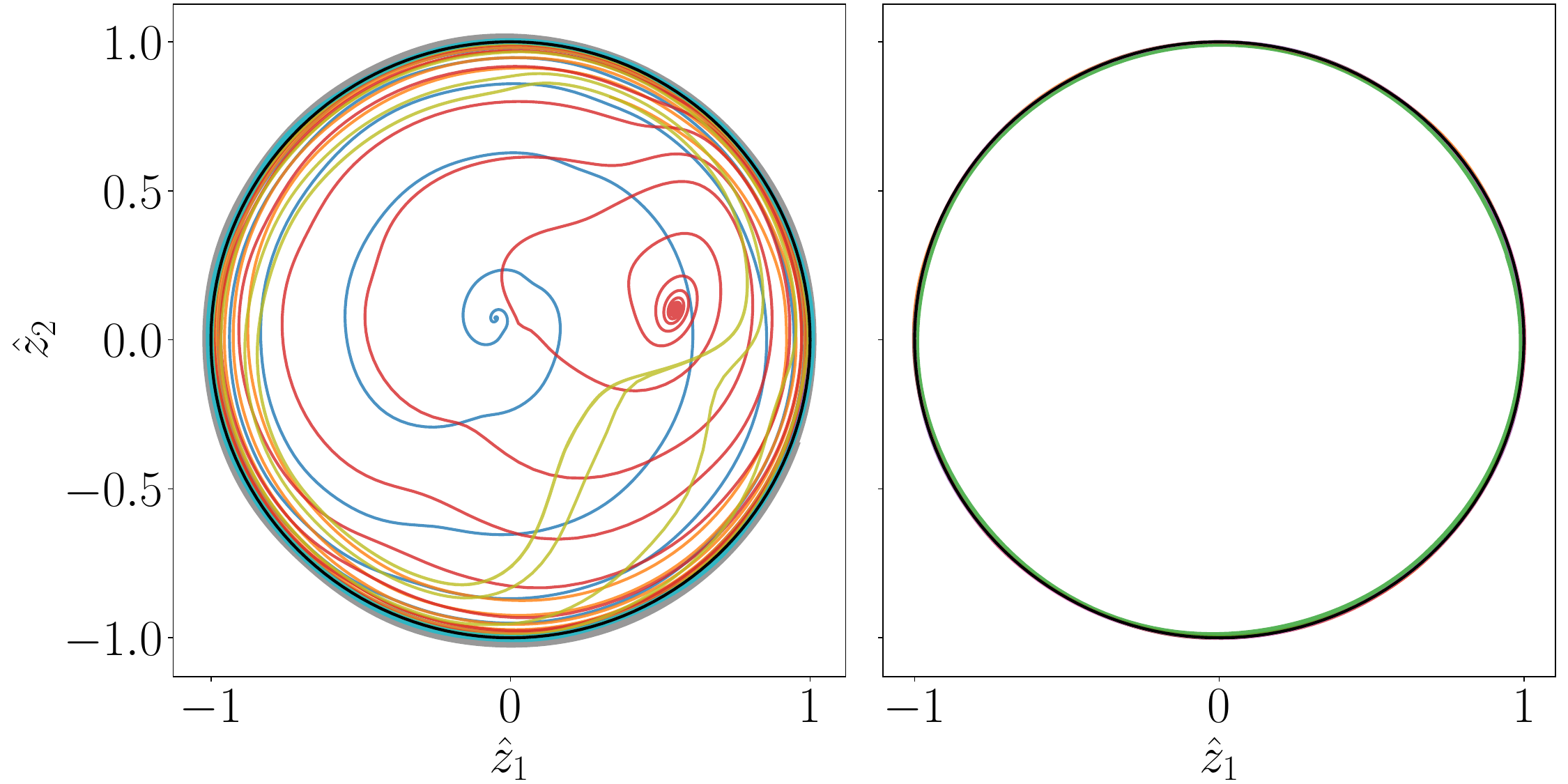}
    
     \caption{Normalized latent trajectory produced by 10 classic ROMs (left) and 10 symmetric ROMs (right)} 
     \label{fig: Circle traj}
\end{figure}
To evaluate the robustness of the ROMs, based on classical and symmetric CAEs, we generated $10$ instances of the classic ROMs and $10$ instances of the symmetric ROMs to assess the reproducibility of the results. With a little abuse of notation we denote with classic and symmetric ROMs, the ROMs that employ classic and symmetric CAEs, respectively. For the symmetric variant, we utilize a parameter-free decoder that serves as the pseudoinverse of the encoder (i.e. $\boldsymbol{f} = 0$ in Equation \ref{eq: projecion layer}), motivated by the simplicity of the problem at hand. 

To facilitate a direct comparison between models trained independently, which may produce topologically equivalent but geometrically distinct latent representations, we map all latent trajectories onto a unit circle. For a high-dimensional trajectory $\boldsymbol{u}_h(t, \boldsymbol{\mu})$, first we project it onto the latent space ($\boldsymbol{z}(t, \boldsymbol{\mu}) = E(\boldsymbol{u}_h(t, \boldsymbol{\mu}))$), and then define the trajectory center of mass $\bar{\boldsymbol{z}}(\boldsymbol{\mu})$ as:
\begin{equation}
\bar{\boldsymbol{z}}(\boldsymbol{\mu}) = \frac{\int_0^T \boldsymbol{z}(t, \boldsymbol{\mu}) \| \dot{\boldsymbol{z}}(t, \boldsymbol{\mu}) \| dt}{\int_0^T \| \dot{\boldsymbol{z}}(t, \boldsymbol{\mu}) \| dt}.
\end{equation}
Next, introduce polar coordinates relative to the center:
\begin{equation*}
    \theta(\boldsymbol{x})=\text{atan2}{\left(\frac{x_2}{x_1}\right)}, \qquad \boldsymbol{x}\in \mathbb{R}^2
\end{equation*}
and define the radial profile of the trajectory as the function
\begin{equation*}
    \rho(\alpha)=\| \boldsymbol{z}(t, \boldsymbol{\mu}) - \bar{\boldsymbol{z}}(\boldsymbol{\mu}) \|, \quad \text{such that}  \quad \theta(\boldsymbol{z}(t, \boldsymbol{\mu}) - \bar{\boldsymbol{z}}(\boldsymbol{\mu}))=\alpha.
\end{equation*}
In other words, $\rho(\alpha)$ gives the distance from the origin of the true latent trajectory along direction $\alpha$.
It is now possible to map each point $\boldsymbol{z}(t, \boldsymbol{\mu})$ into a unique point of the unit circle $\hat{\boldsymbol{z}}$, by means of a normalization function $\mathcal{N}$:
\begin{equation}
\hat{\boldsymbol{z}}(t, \boldsymbol{\mu}) = \frac{\boldsymbol{z}(t, \boldsymbol{\mu}) - \bar{\boldsymbol{z}}(\boldsymbol{\mu})}{\rho(\theta(\boldsymbol{z}(t, \boldsymbol{\mu}) - \bar{\boldsymbol{z}}(\boldsymbol{\mu})))} = \mathcal{N}(\boldsymbol{z}(t, \boldsymbol{\mu}); E, \boldsymbol{u}_h). 
\end{equation}
This normalization function $\mathcal{N}$ is solely determined by the encoder $E$ and the reference solution $\boldsymbol{u}_h$, as the shifting constant $\bar{\boldsymbol{z}}(\boldsymbol{\mu})$ and the scaling function $\rho$ are defined by the true latent trajectory $\boldsymbol{z}(t, \boldsymbol{\mu})$. More generally, for any point $\boldsymbol{x}\in \mathbb{R}^2$ we define
\begin{equation*}
    \boldsymbol{\hat{x}} =  \mathcal{N}(\boldsymbol{x}; E, \boldsymbol{u}_h) = \frac{\boldsymbol{x} - \bar{\boldsymbol{z}}(\boldsymbol{\mu})}{\rho(\theta(\boldsymbol{x} - \bar{\boldsymbol{z}}(\boldsymbol{\mu})))}.
\end{equation*} 
Notably, the normalization function maps the predicted ROM trajectories onto the unit circle only when they accurately follow the ground truth. Any deviation from the true latent trajectory is thus reflected as a departure from the unit circle of the normalized ROMs prediction.

After the training phase, the ROMs were evolved for $10$ periods ($t \in [0, 10T]$) to test the stability of the latent flow. 
\begin{figure}[h!]
    \centering
    \includegraphics[width=\linewidth]{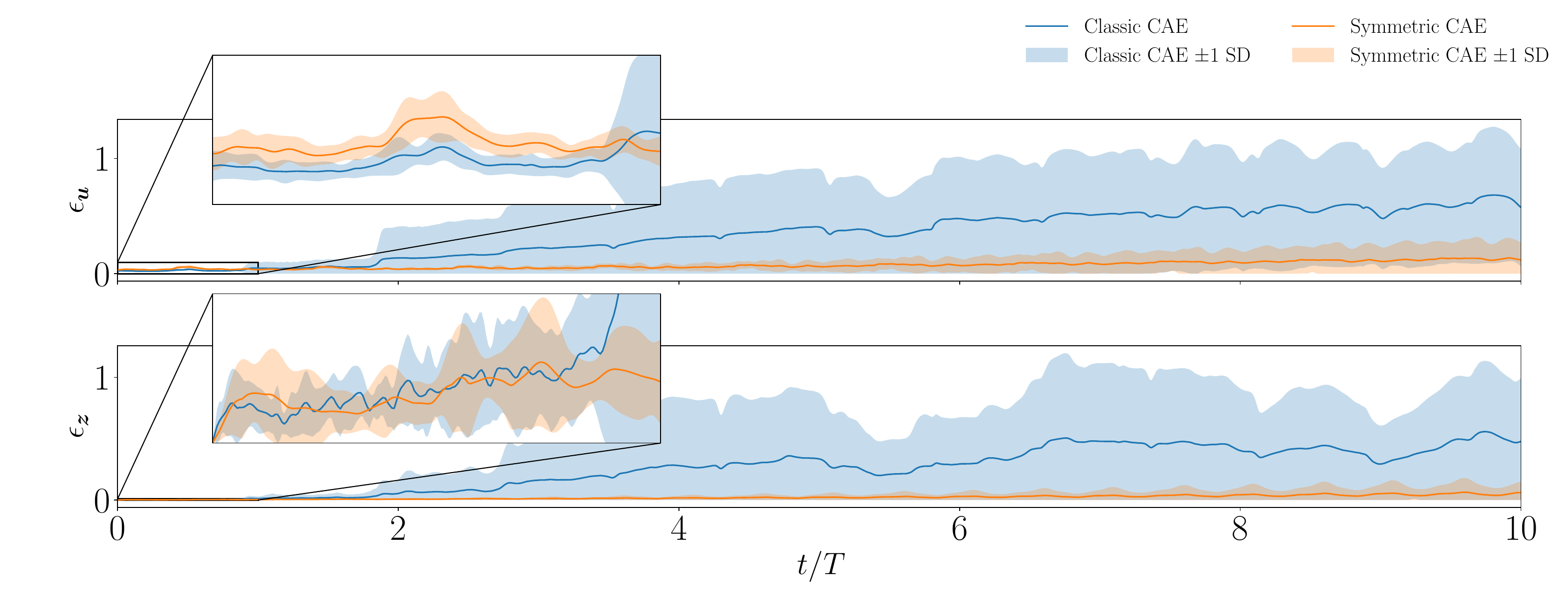}
    \caption{Statistical evolution of the high-dimensional reconstruction error (top) and latent space error (bottom) for the classic and symmetric ROMs. Results are computed over a 10-period time evolution across the test set. Solid lines represent the mean error, while the shaded regions indicate the standard deviation.}
    \label{fig:LA evolution error}
\end{figure}
The resulting normalized trajectories are visualized in Figure \ref{fig: Circle traj}. The $10$ trajectories generated by the symmetric ROMs tend to stay closer to the unit circle, where they should belong. This is not the case for the classic one because, soon after the second period, some of the trajectories depart from the unit circle and start wandering around. A clear evolution of the error defined by Equation \ref{eq: epsilon u} and the error in latent space, defined as
\begin{equation}
    \epsilon _ {\boldsymbol{z}}(t) = \| \hat{\boldsymbol{z}}(t) - \mathcal{N}(E(\boldsymbol{u}_h(t))) \|
\end{equation}
is illustrated in Figure \ref{fig:LA evolution error}. The symmetric ROMs maintain an error profile that is significantly lower than that of the classical models. We attribute this superior performance to three primary factors:
\begin{enumerate}
    \item Having a pair of encoder-decoder models that are consistent provides a latent space that is a true coordinate system. Any latent point which is decoded from the latent space to the physical space and then encoded it again, it results in the same input point.
    \item The constraint $E \circ D = id$ acts as a powerful regularizer during training. This leads to a smoother latent manifold that is inherently more amenable to interpolation and temporal integration by the latent flow network.
    \item The symmetric constraint promotes a lower Lipschitz constant for the latent flow \cite{farenga2025latent}. This reduces the sensitivity of the dynamics to small numerical perturbations, ensuring long-term stability in periodic regimes.
\end{enumerate}

\subsection{Viscous Burger Equation}
In the second test case, we explore the application of symmetric CAEs in the context of non-linear parametric PDE. Given that the only two parameters involved are time and viscosity, we can choose a two-dimensional latent space size ($l=2$). Multiple models were trained, but we focus on the two best-performing models. 

\begin{figure}[h!]
    \centering
    \begin{subfigure}{0.4\linewidth}
        \centering
        \includegraphics[width=1\linewidth]{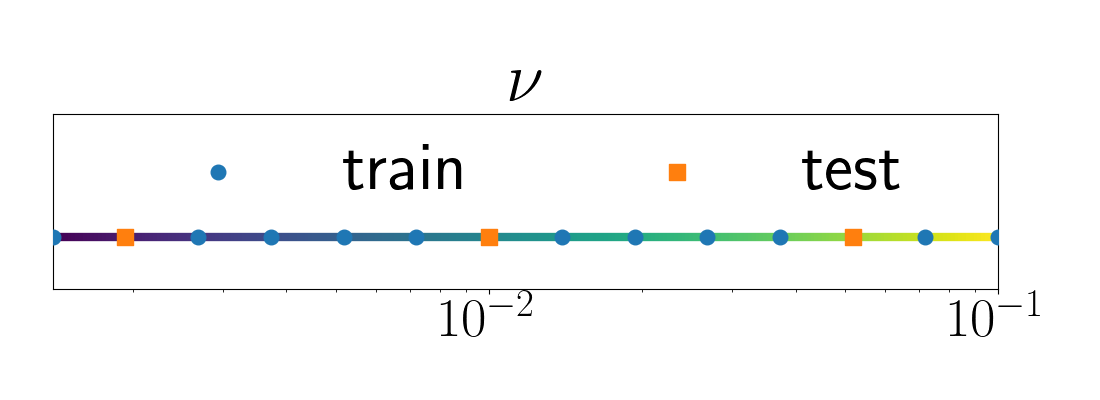}      
    \end{subfigure} \\
    \centering
    \begin{subfigure}{0.45\linewidth}
        \centering
        \includegraphics[width=1\linewidth]{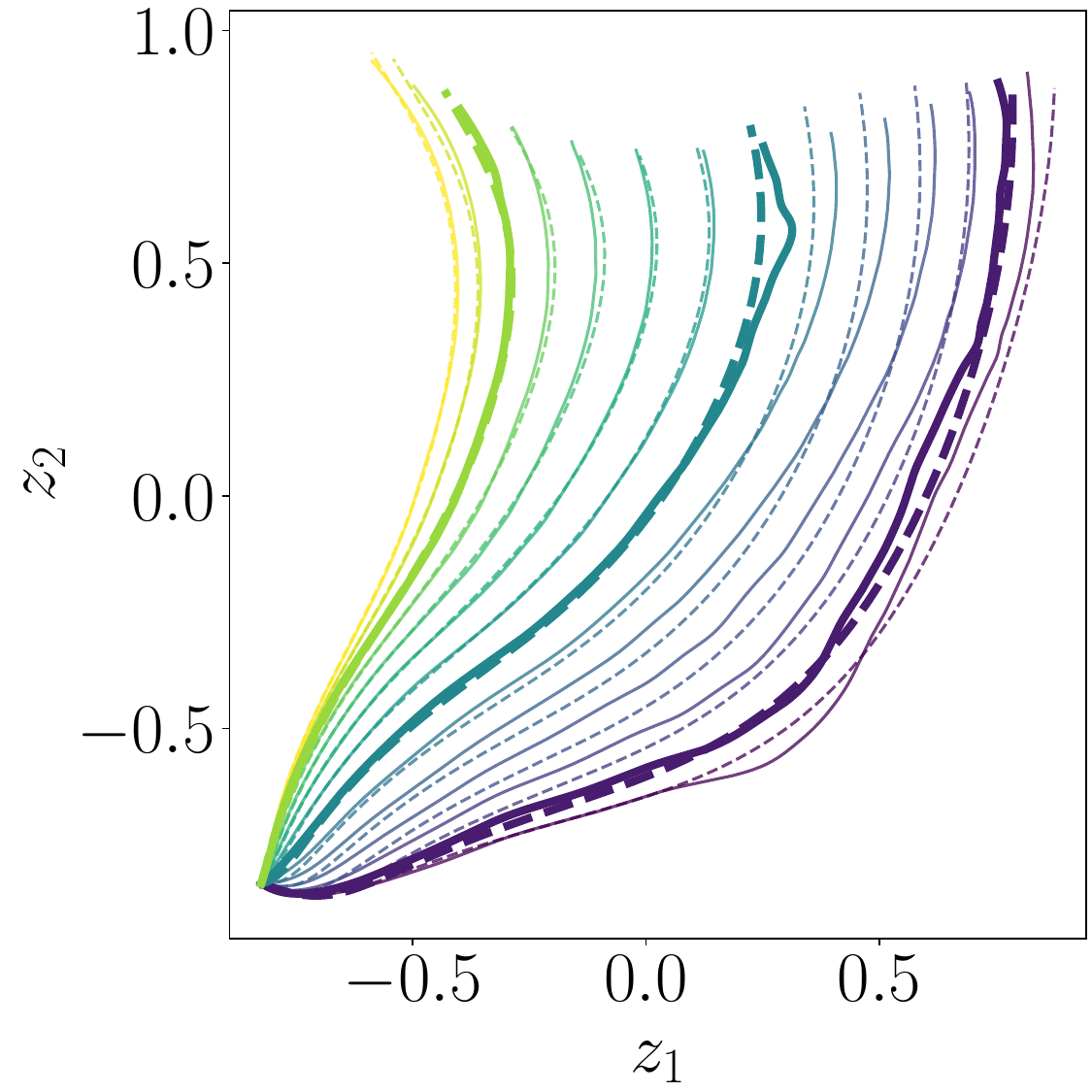}
        \caption{}
    \end{subfigure}
    \begin{subfigure}{0.45\linewidth}
        \centering
        \includegraphics[width=1\linewidth]{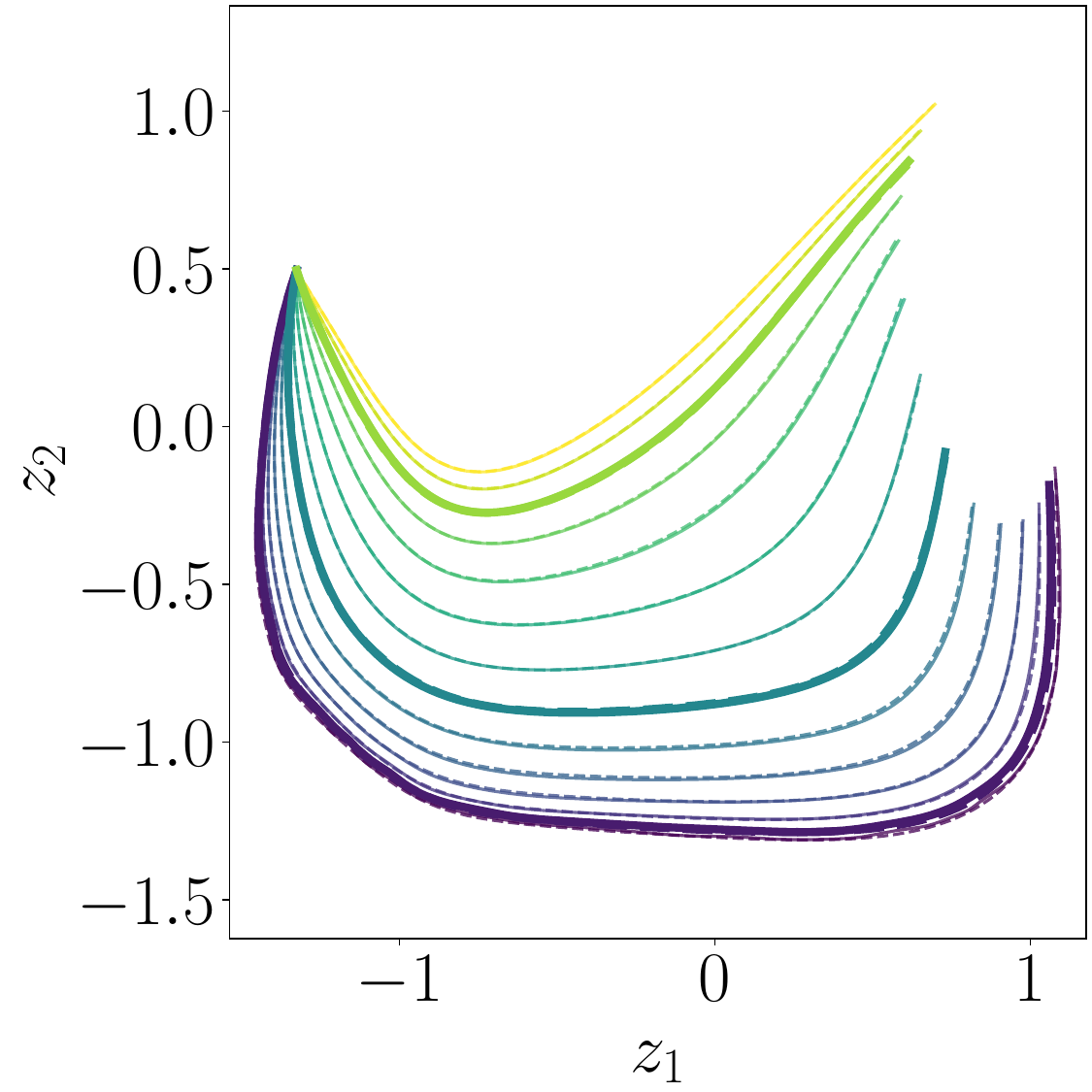}        
        \caption{}
    \end{subfigure}

    \caption{Comparison of the latent space trajectories ($l=2$) for the parametric Viscous Burgers' equation using Classic CAE (left) and Symmetric CAE (right). Solid lines represent the projection of the high-fidelity solution through the encoder, $\mathbf{z} = E(\mathbf{u}_h)$, while dashed lines denote the trajectories predicted by the ROM. The trajectories are colored according to the viscosity $\nu$, with the thicker lines highlighting the three cases from the unseen test set.}
    \label{fig: VB LS}
\end{figure}

Figure \ref{fig: VB LS} illustrates the reduced trajectories in the latent space. The lack of smoothness in the latent space produced by the classical CAE is particularly evident when examining viscosity values from the test trajectories. The solid lines in Figure \ref{fig: VB LS} represent the projection of the FOM solution \( \boldsymbol{z} = E(\boldsymbol{u}_h) \) through the encoder. Among these, three lines stand out as notably thicker: these correspond to the test set. When employing a standard convolutional AE for compression, the test trajectories exhibit significant wiggling, indicating irregularities in regions distant from the training points. In contrast, when projecting the curve using the Symmetric CAE, these irregularities are substantially reduced, further demonstrating the regularizing effect inherently introduced by this architecture.

This issue also impacts the chain of prediction capabilities within the ROM. Specifically, the LFNet relies on the knowledge of the tangent vector (the derivative taken along the trajectory) at a particular point in the latent space. When the derivative of a given noisy curve is computed, it tends to amplify the effects of the noise, resulting in even more irregular derivatives (tangent vectors). These perturbations cause the predicted trajectory by classic ROM (computed via Equation \ref{eq: ROM solution}) and depicted by the dashed lines in the left plot of Figure \ref{fig: VB LS}, to deviate from the reference trajectory due to the accumulation of large errors. In this framework, the overlap of the solid and dashed lines for the symmetric ROM can be explained by better accuracy in approximating the smoother time derivative, leading to slower error propagation. 

\begin{figure}[h!]
    \centering
    \includegraphics[width=0.75\linewidth]{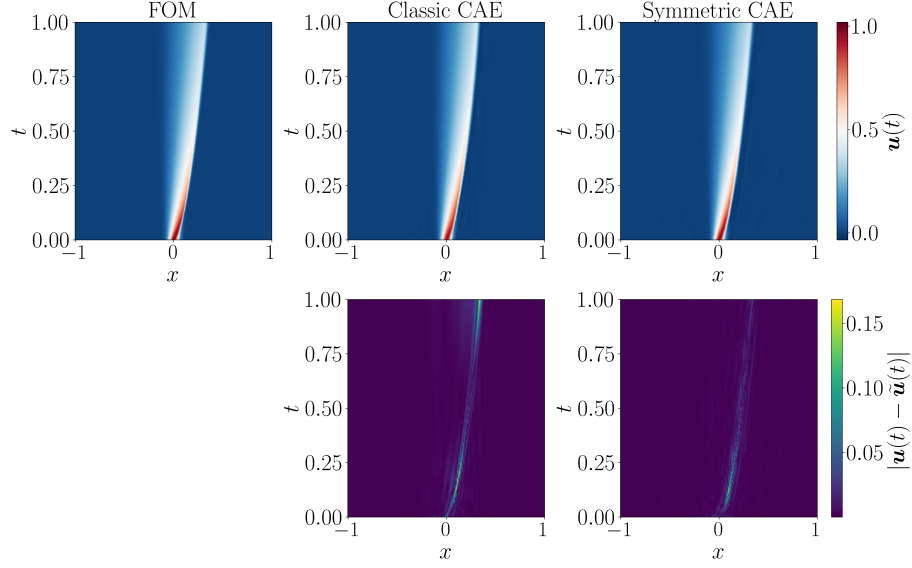}
    \caption{First line: spatio-temporal solutions generated by the FOM (left), the classical ROM (center), and the symmetric ROM (right), using a viscosity value $\nu$ selected from the training set. Second line: point-wise error as functions of the spatial and temporal coordinates for classical ROM (center) and symmetric ROM (right).}
    \label{fig: VB training}
\end{figure}

\begin{figure}[h!]
    \centering
    \includegraphics[width=0.75\linewidth]{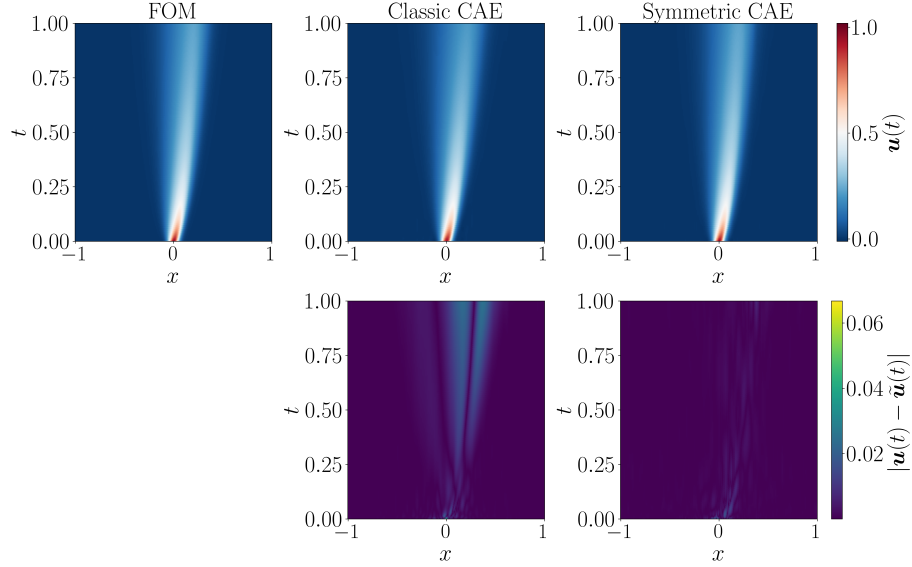}
    \caption{First line: spatio-temporal solutions generated by the FOM (left), the classical ROM (center), and the symmetric ROM (right), using a viscosity value $\nu$ selected from the test set. Second line: point-wise error as functions of the spatial and temporal coordinates for classical ROM (center) and symmetric ROM (right).}
    \label{fig: VB test}
\end{figure}

\begin{figure}[h!]
    \centering
    \includegraphics[width=1\linewidth]{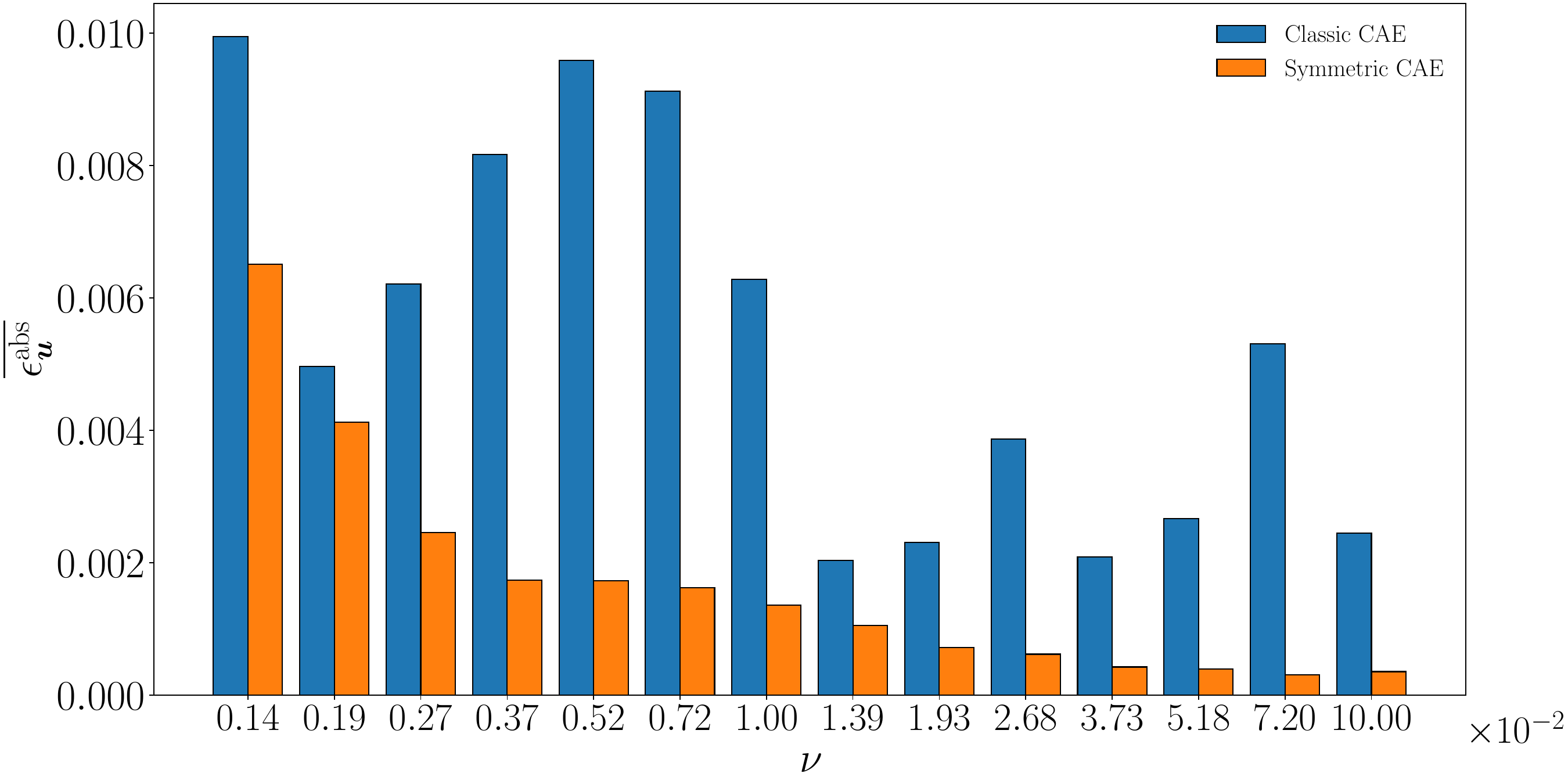}
    \caption{Error $\overline{\epsilon_{\boldsymbol{u}} ^{\text{abs}}} $ defined as the time mean of Equation \ref{eq: absolute rec error}, for ROM based on standard CAE (left) and symmetric CAE (right).}
    \label{fig: VB error}
\end{figure}

The improved regularity of the latent trajectories translates into superior physical reconstruction. Figures \ref{fig: VB training} and \ref{fig: VB test} compare the FOM solutions with the ROM predictions and their corresponding punctual errors for viscosity values in both the training and testing sets. For the chosen training viscosity value (Figure \ref{fig: VB training}), the classical architecture notably fails to capture the precise position of the shock front, the most fundamental feature of the flow in this regime. Furthermore, as the viscosity increases and the importance of the non linear term decreases accordingly, the diffusive nature of the PDE results in a more smeared shock profile. In these higher-viscosity test cases (Figure \ref{fig: VB test}), the classical ROM manifests significantly larger errors for $t \ge 0.2$, precisely when the diffusion process begins to dominate. The symmetric ROM, however, maintains high fidelity across both regimes.
In order to make a general comparison, we introduce the absolute error:
\begin{equation}
    \label{eq: absolute rec error}
    \epsilon_{\boldsymbol{u}} ^{\text{abs}} (t; \nu) = \|\boldsymbol{u}_h(t; \nu) - D(\boldsymbol{\hat{z}}(t; \nu)) \|_{L_2}.
\end{equation}
The performance of the classical and symmetric ROMs is evaluated via the mean absolute error, defined as:
\begin{equation*}
\overline{\epsilon_{\boldsymbol{u}}^{\text{abs}}}(\nu) = \frac{1}{T} \int_0^T \epsilon_{\boldsymbol{u}}^{\text{abs}}(t; \nu) , dt.
\end{equation*}
As illustrated in Figure \ref{fig: VB error}, a distinct trend characterizes the symmetric architecture: the mean error consistently decreases as the viscosity parameter increases across both training and testing ranges. By contrast, the classical ROM lacks this monotonic behavior, displaying a less predictable error distribution across the parametric space.
Detailed analysis of the temporal evolution reveals the main source of these errors. In low-viscosity regimes, the classical ROM fails to resolve precise shock locations, resulting in peak error values at the beginning of the simulation. This behavior is clearly visible in the bottom row of Figure \ref{fig: VB training}. Conversely, in viscosity-dominated scenarios, the classical model exhibits its largest inaccuracies toward the end of the temporal domain, as shown in Figure \ref{fig: VB test}. While the symmetric ROM follows similar qualitative patterns, it maintains significantly lower error magnitudes on average. This superior performance, summarized in Figure \ref{fig: VB error}, underscores the greater reliability preserved by the symmetric CAE-ROM throughout the investigated parametric range.

\subsection{Kuramoto-Sivashinsky Equation}

The application of ROMs to chaotic systems, such as the Kuramoto-Sivashinsky equation, represents a significant challenge due to the inherent unpredictability of the dynamics. In chaotic regimes, trajectories exhibit a sensitive dependence on initial conditions, where small perturbations grow exponentially. Consequently, instantaneous field comparisons, standard in non-chaotic cases, lose their physical significance after a short integration interval. This unpredictability adds complication in understanding weather the reduced model is able to actually capture the true essence of the system being modeled (invariant properties of the chaotic attractor). We should then be equipped with tools that reveal the essence of the dynamics and extract quantities that are truly representative of the chaotic system at hand. These tools can be recast into two macro groups: stability analysis (Covariant Lyapunov vectors and Lyapunov exponents) and statistical analysis (evolution of the mean variables, Probability Density Functions (PDFs), energy spectrum, etc.). The latter are extensively employed to analyze turbulent fluid fields, which can be considered one of the prime examples of chaos \cite{pope2001turbulent}.

In this study, we fix the KS parameters to $L = 20\pi$ and $\nu = 1$ to ensure a fully developed chaotic regime \cite{ozalp2025stability}, characterized by a leading Lyapunov Exponent (LE) $\lambda = 0.085$, which is associated to a characteristic Lyapunov Time (LT) $\tau_ \lambda = 1/\lambda$. The space domain is discretized using $N_h = 512$ nodes, and the equation is solved using a time step $\Delta t = 0.25$, for a total time of $T= 2.5 \times 10^4$, corresponding to $2125$ LTs. We chose a latent space size $l=24$, based on the study made by Özalp and Magri \cite{ozalp2025stability}. To evaluate the performance of the classic and symmetric ROMs, we move beyond trajectory-based errors and assess the models' ability to capture the flow's statistical and spectral characteristics. 

The classical CAEs exhibited a more pronounced lack of robustness when subjected to the chaotic regime of the KS equation (compared, for instance, to the Linear Advection test case). Of the three instances trained, two suffered from instability, with their solutions effectively ``exploding'' after only a few LTs. In these instances, the norm of the latent state solution $\boldsymbol{z}$—calculated via Equation \ref{eq: ROM solution}—increases rapidly, escaping the latent attractor (either the global attractor or the inertial manifold) showing an unphysical behavior and drifting into regions of the latent space that were never encountered during the training phase. Decoding these regions far from the inertial manifold, leads to an unphysical high-dimensional solution that depends strongly on the weights learned and it is unpredictable in principle. Even the single instance that achieved partial stability remained highly sensitive, manifesting instabilities for several specific initial conditions.

In contrast, the models based on the symmetric AE architecture proved to be remarkably robust. All seven 
generated instances maintained stable behavior throughout the integration period, representing a significant improvement over the classical approach. Therefore, the consistent stability of the symmetric models confirms that the structural constraints effectively regularize the latent flow, keeping the dynamics confined to the attractor.

The prediction horizon is defined as the temporal interval beyond which the ROM trajectory noticeably diverges from the ground-truth high-fidelity solution. In the existing literature, surrogate models for the KS equation—particularly those utilizing autoregressive approaches on FOM snapshots (e.g., RNNs, LSTMs, and Reservoir Computing)—report prediction horizons ranging from $1$ to $5$ Lyapunov Times \cite{vlachas2020backpropagation, pathak2018model, pathak2018hybrid}. When Reservoir Computing is applied specifically to latent dynamics, this horizon is often reported around $1.25$ LT \cite{pathak2018model}.
\begin{figure}[h!]
    \centering
    \includegraphics[width=0.75\linewidth]{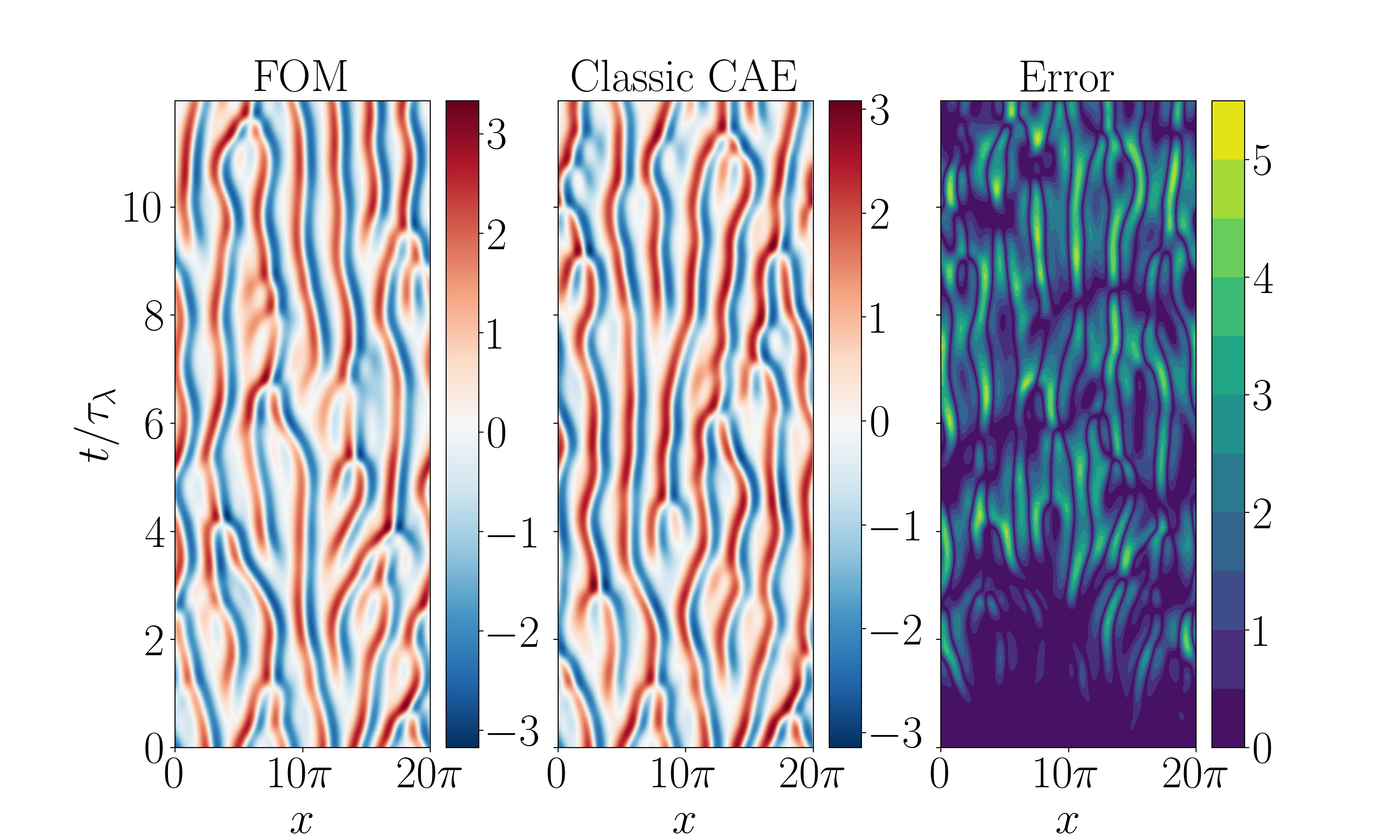}
    \caption{Reference solution of the KS equation (left), prediction of the ROM model based on classical CAE (center) and point-wise absolute error (right).}
    \label{fig: KS_fieldCls}
\end{figure}
\begin{figure}[h!]
    \centering
    \includegraphics[width=0.75\linewidth]{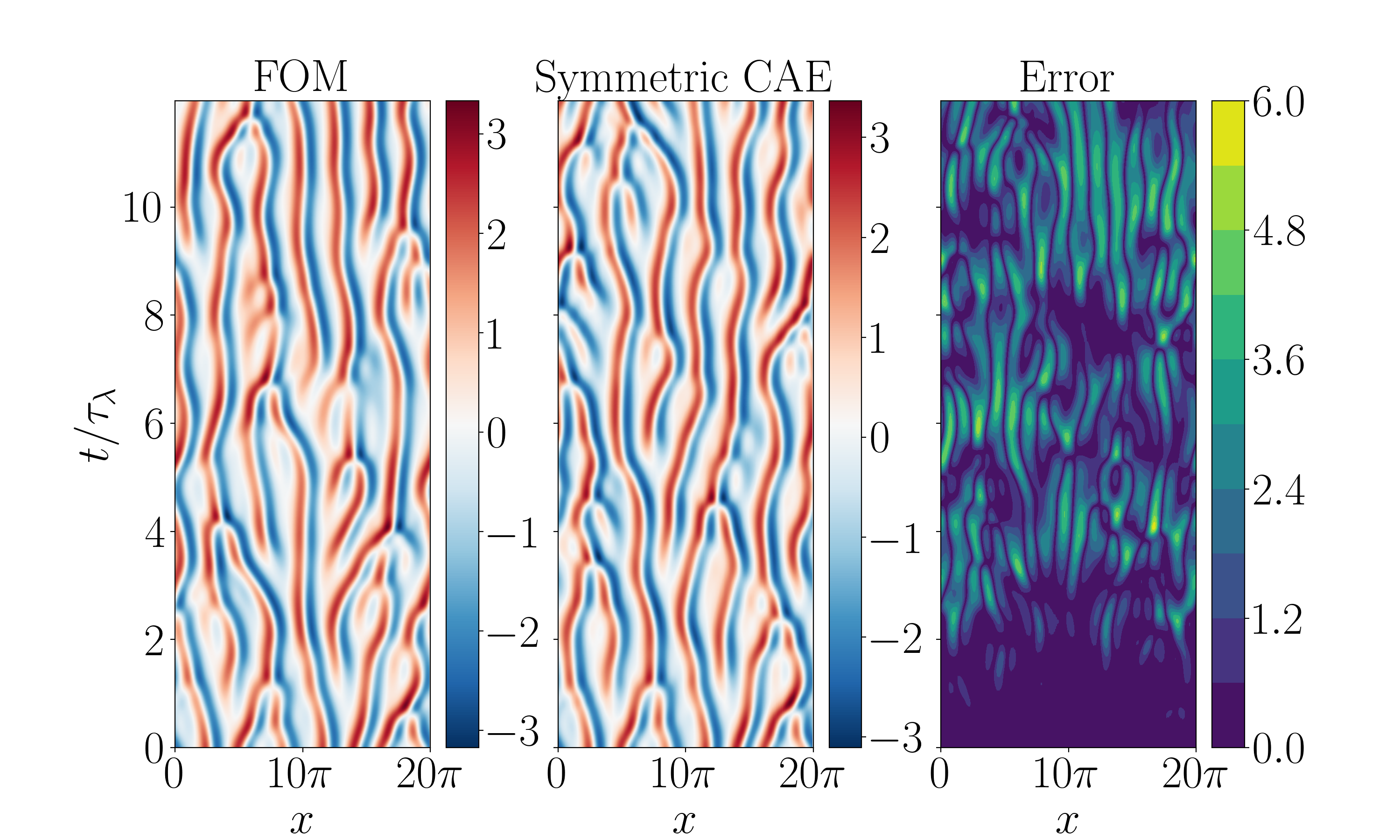}
    \caption{Reference solution of the KS equation (left), prediction of the ROM model based on symmetric CAE (center) and point-wise absolute error (right).}
    \label{fig: KS_fieldSym}
\end{figure}
Starting from a selected initial condition, we integrated the reduced-order models—utilizing both the classical CAE (Figure \ref{fig: KS_fieldCls}) and the symmetric architecture (Figure \ref{fig: KS_fieldSym})—for a total duration of $12$ LT. Both models display the characteristic features of the KS solutions, with strides that merge and split over time. The prediction horizon stands at about $1$ and $2$ LT for classic and symmetric ROMs respectively, which is not too far from $1.25$ found by Özalp and Magri \cite{ozalp2025stability}. 

The decision to employ a Latent Flow Network as the temporal propagator, rather than traditional Recurrent Neural Networks or Reservoir Computing, provides distinct advantages in terms of both training efficiency and model robustness. Primarily, the LFNet avoids the need for Backpropagation Through Time, which is a requirement for training LSTM or GRU-like architectures. While increasing the unroll length in RNNs can enhance predictive performance, it simultaneously increases memory overhead and computational cost, often rendering long-range dependency training prohibitive \cite{werbos2002backpropagation}. 
Furthermore, while Reservoir Computing methods circumvent traditional problems of BPTT \cite{vlachas2020backpropagation}, their performance is notoriously sensitive to the calibration of hyperparameters, such as the spectral radius and Tikhonov regularization parameter \cite{margazoglou2023stability}. By shifting the objective to the supervised learning of the latent state's time derivative, the LFNet is a ``plug-and-play'' surrogate. This approach allows the network to internalize the underlying dynamics during the training phase, effectively eliminating the exhaustive parameter tuning associated with RC while maintaining equal or even superior predictive accuracy.
\begin{figure}[h!]
    \centering
    \includegraphics[width=0.75\linewidth]{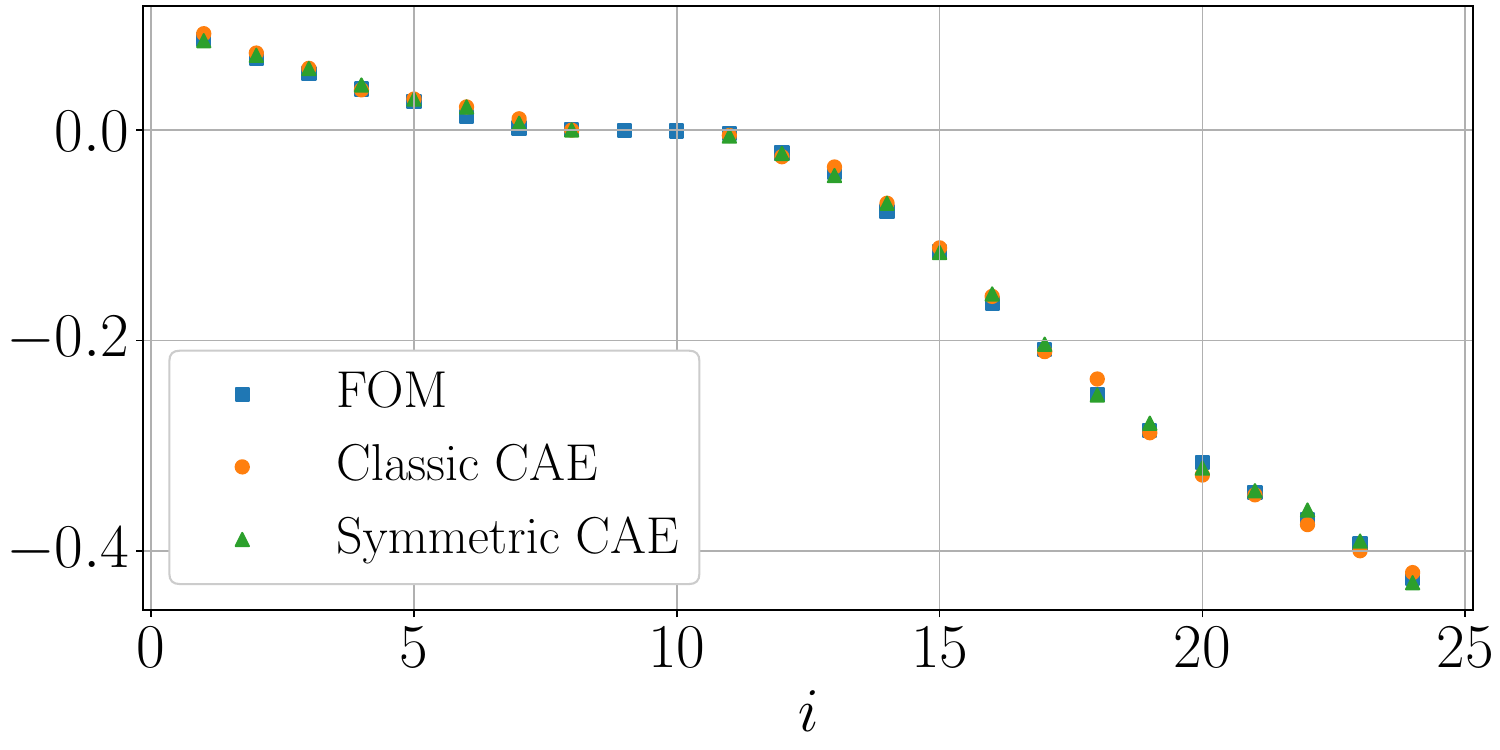}
    \caption{Lyapunov Spectrum predicted by the classical and symmetric ROMs and the reference spectrum obtained from the high fidelity simulation.}
    \label{fig: LEs}
\end{figure}

Another useful analysis concerns the Lyapunov exponents, , which measure the sensitivity of a solution to initial conditions by quantifying the average rate of growth or decay of a small perturbation along a given direction. For a rigorous theoretical treatment and numerical algorithms to compute LEs, we refer the reader to \cite{skokos2009lyapunov, oseledec1968multiplicative, benettin1980lyapunov}. An high level overview about the methods to compute the LEs from the FOM equation is provided in Appendix \ref{sec: FOM LEs}.
To compute the LEs of the ROM, we leverage the automatic differentiation capabilities of PyTorch \cite{paszke2017automatic}. Given that the latent dynamics of the ROM are governed by 
\begin{equation}
\label{eq: ROM eq LE}
    \dot{\mathbf{z}}(t; \boldsymbol{\mu}) = h_\theta(\mathbf{z}, \boldsymbol{\mu}, t) 
\end{equation}
the perturbation evolves according to:  
\begin{equation}
\label{eq: ROM variational equation}
    \dot{\mathbf{v}}(t) = d_ {\boldsymbol{z}} h_\theta(\mathbf{z}, \boldsymbol{\mu}, t) \mathbf{v}
\end{equation}
where $d_ {\boldsymbol{z}} h_\theta(\mathbf{z}, \boldsymbol{\mu}, t)$ can be easily computed using backpropagation. Equation \ref{eq: ROM eq LE} and \ref{eq: ROM variational equation} are solved together using a RK4 scheme, with time stepping for equation \ref{eq: ROM eq LE} which is half of the time step used to solve equation \ref{eq: ROM variational equation}.
Furthermore, Appendix \ref{sec: Equivalence LE} provides a mathematical proof demonstrating that, under specific conditions, the LEs of the ROM and FOM are equivalent. As established in Appendix \ref{sec: Lyapunov Exponents}, the LS serves as the ``fingerprint'' of a chaotic system, encapsulating critical information regarding its dimensionality, entropy, and long-term predictability. Figure \ref{fig: LEs} illustrates the comparison of the LEs computed by the FOM with those derived from the ROMs. The models accurately capture the positive LEs. However, as noted in \cite{pathak2017using, vlachas2020backpropagation}, ROMs are unable to represent the multiplicity of the zero LEs. This behavior is suggested in \cite{pathak2017using} but has not been elaborated upon. As discussed in Appendix \ref{sec: Symmetries and LEs}, the zero LEs are linked to the number of continuous symmetries present in the dynamical system. In particular, the Kuramoto-Sivashinsky equation exhibits three continuous symmetries: time-translation invariance, space-translation invariance, and Galilean invariance. Neural networks do not maintain symmetries such as space translation or Galilean invariance \cite{weiler2023EquivariantAndCoordinateIndependentCNNs}. Consequently, the only symmetry retained in the latent dynamical system is the time-translation invariance associated with the Lyapunov vector $\dot{\mathbf{z}}(\mathbf{z},t ; \boldsymbol{\mu})=d\boldsymbol{x}_\alpha ^{-1}(\mathbf{f}_h\left(t, \mathbf{x}_\alpha (\mathbf{z}(t ; \boldsymbol{\mu})) ; \boldsymbol{\mu}\right))$. Given that the ROMs compute $l$ LEs while omitting the two zero LEs, it is evident that the ROM produced two spurious LEs. By manually eliminating them, you can observe the overlapping between the true and predicted spectra in Figure \ref{fig: LEs}.

\begin{figure}[h!]
    \centering
    \includegraphics[width=0.75\linewidth]{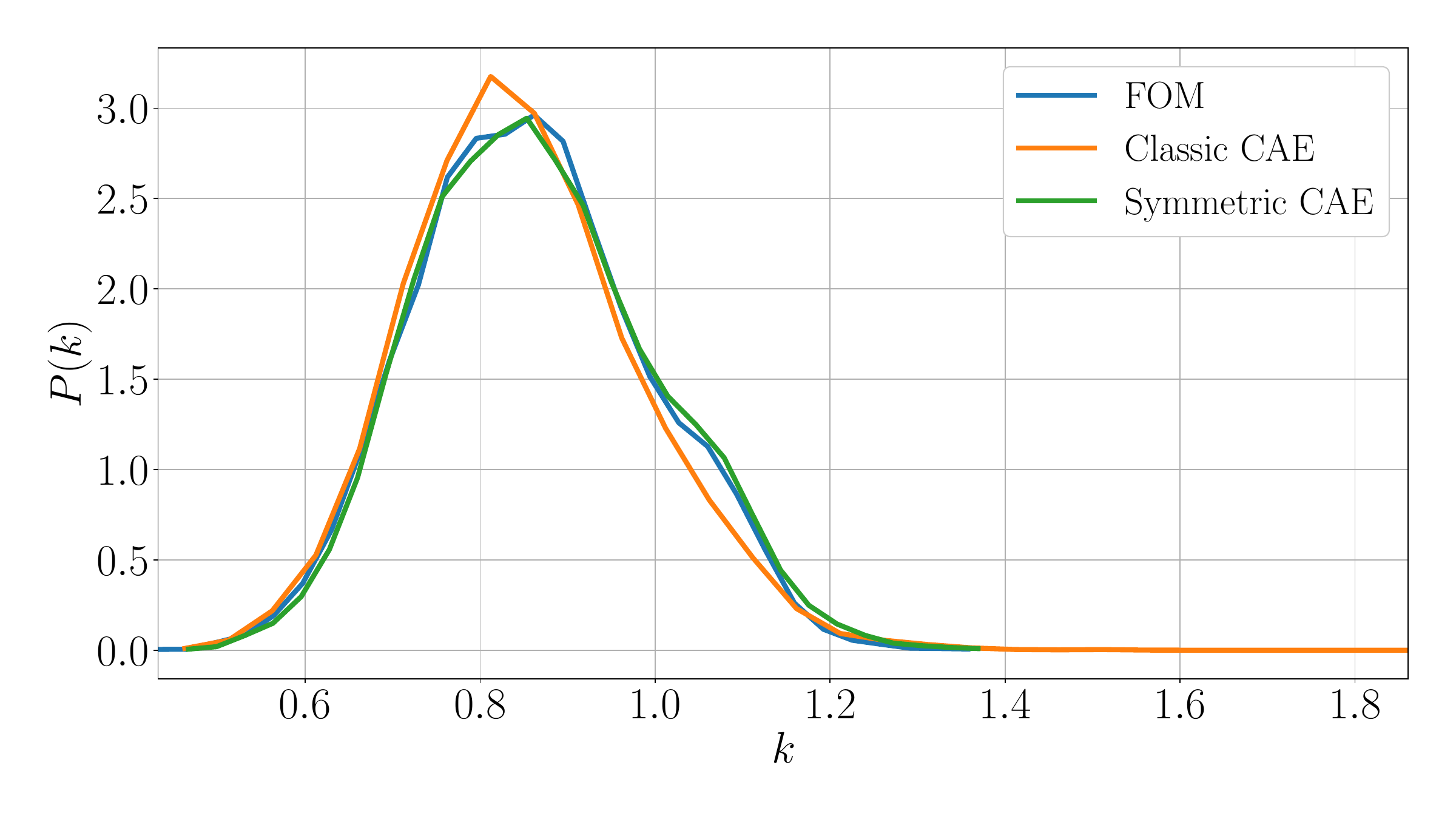}
    \caption{Probability Density Function (PDF) of the kinetic energy.}
    \label{fig: k PDF}
\end{figure}

We computed the Kaplan-Yorke (KY) dimension taking into account the two missing LEs. The KY dimension and the value of the first LE is reported in table \ref{tab: KY and LE}, along with the relative errors.
\begin{table}[]
    \centering
    \begin{tabular}{|c|c|c|c|}
    \hline
        & FOM & Classic & Symmetric  \\
        \hline
     KYD    & $16.18$ & $16.48$ & $16.37$\\
     \hline
     rel. error [\%] &$ -$ & $1.87$ & $1.84$  \\
     \hline
     $\lambda_1$ & $0.08567$ & $0.09153$ & $0.08489$ \\
     \hline
     rel. error [\%] & $-$  &$6.83$ & $0.93$ \\
     \hline
    \end{tabular}
    \caption{Comparison of FOM and ROMs leading LE and KY dimension.}
    \label{tab: KY and LE}
\end{table}
To assess the long-term ergodic properties of the reduced-order models, we computed the Probability Density Function of the kinetic energy over an extensive integration period of $2083$ LT. The kinetic energy is defined as:
\begin{equation*}
    k(t) = \frac{1}{2} \int _ 0 ^  L  u^2(x, t) dx.
\end{equation*}
The statistical characterization of the classical ROM was notably hampered by its inherent numerical instability; as previously discussed, certain initial conditions triggered divergent trajectories. These unstable instances resulted in unphysical state reconstructions, which produced an erroneous PDF with a pronounced bias toward high-energy states. As illustrated in Figure \ref{fig: k PDF}, while both the classical and symmetric ROMs capture the general distribution of the high-fidelity kinetic energy, the classical model exhibits spurious high-energy configurations that are absent in the FOM data.
To conclude the spectral characterization of the Kuramoto-Sivashinsky system, we compute the averaged kinetic energy spectrum utilizing Welch’s method. This analysis allows for a precise evaluation of how the reduced-order models distribute energy across different spatial scales. The symmetric ROM demonstrates a superior capability in capturing the energy content across the entire frequency spectrum, whereas the classical CAE-based ROM exhibits significant discrepancies. These inaccuracies are particularly evident in the low-frequency regime, which corresponds to the large-scale coherent structures.
\begin{figure}[h!]
    \centering
    \includegraphics[width=0.75\linewidth]{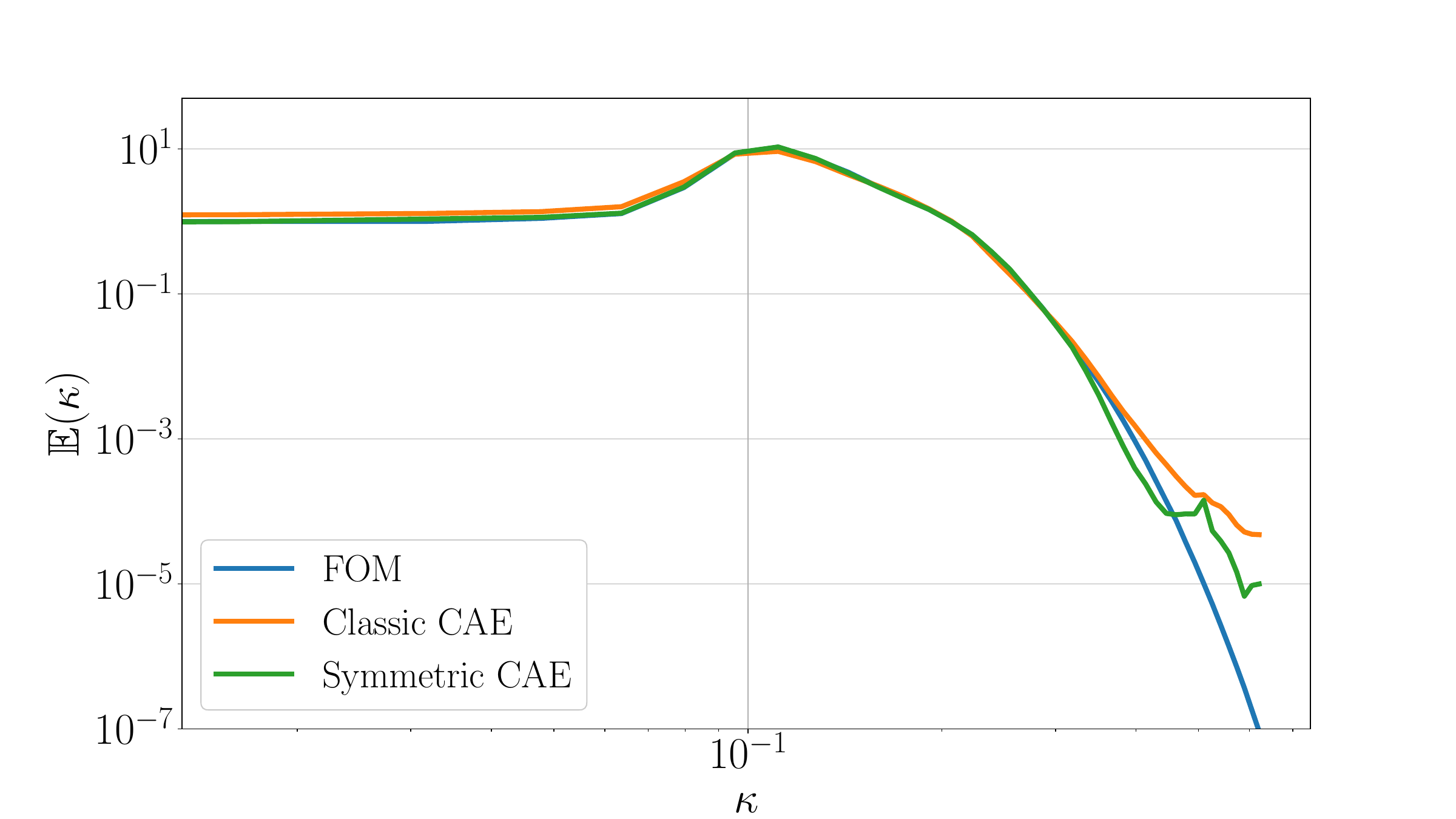}
    \caption{Kinetic energy spectrum against wave number for the FOM, classical and symmetric ROM systems.}
    \label{fig: PSD}
\end{figure}
%




\section{Conclusions}
\label{sec: Conclusion}
In this work, we begin by mathematically deriving the equation of the reduced model, leveraging some basic knowledge from differential geometry. Next, we adopt tools from the machine learning domain that, without any external assistance—i.e., in an unsupervised manner—are capable of learning an approximation of the manifold's parametrization by examining a subset of points belonging to that manifold. Unfortunately, the standard versions of these tools do not fully satisfy all the properties that parametrization should possess, particularly the concept of representation consistency. Recent literature has introduced augmented neural network types, specifically autoencoders with fully connected layers, that address the representation consistency property. This paper extends this property to hybrid convolutional autoencoders. In this context, ``hybrid'' refers to an encoder composed of convolutional layers, while the decoder is constructed with fully connected neural networks designed to be parametrized by a limited number of weights—either matching the encoder or, potentially, none. This approach allows the encoding process to retain the advantages of CNN layers, such as a number of trainable parameters that is independent of the input vector size, shared weights, and a limited receptive field, while permitting the decoder to leverage the full expressiveness of feedforward neural networks for the reconstruction task. Although there are convolutional neural networks capable of processing not only vectors but also 2D arrays, our implementation focuses solely on the 1D version. Addressing the extension to 2D cases presents challenges that will be explored in future work. Consequently, this study is limited to 1D PDEs.

Here, the symmetric AEs, are integrated inside a ROM framework, to see if benefits in terms of performance can be achieved once you get closer to "ideal" maps. The ROM comprises both an AE and a network, which we called Latent Flow Network, able to predict the latent flow at any point and for any given value of the parameters.  Indeed, three different cases were investigated, namely the Linear Advection (LA), the Viscous Burger (VB) equation and the Kuramoto-Sivashinsky (KS) equation. Each of them, carries its own complexity, like advection dominating problem (LA), parameters dependency (VB) or chaos (KS). Throughout these test cases, the symmetric CAE performed equal or better compared to the classical version of the CAEs. In particular, it can be seen that introducing this representation consistency constraint, and in particular the formulation based on projector, acts as a powerful regularizer, producing a latent space which is smoother and easier to be learnt, by any latent dynamic model. Specifically, for LA and VB, ROM based on symmetric AE appeared to be more accurate. Moreover, employing symmetric AEs, results in more stable ROMs (like models that are easier to train or models that doesn't explode). This was particularly clear, when a reduced model was created for the chaotic KS equation. Indeed, it was not so easy to deliver a working model with the classical CAE, while, with the introduction of symmetric AEs we were able to provide viable reduced models for each trained model. For the chaotic KS problem, particular attention is given to the spectra, specifically, to the Lyapunov spectrum. Here we proved first that the FOM and the ROM, possess the same non zeros LEs, and then we explained why, with NN it is difficult to obtain the exact multiplicity of the zero LEs.





\bibliographystyle{abbrv}
\bibliography{sample}
\appendix
\section{Chaotic systems}
\label{sec: chaotic systems}



\subsection{Lyapunov Exponents (LEs) and related quantities}
\label{sec: Lyapunov Exponents}
Let $\boldsymbol{u}(t, \boldsymbol{\mu})$ be the solution of the autonomous dynamical system $\dot{\boldsymbol{u}}(t, \boldsymbol{\mu}) = \boldsymbol{f}(\boldsymbol{u}(t, \boldsymbol{\mu}))$, $\boldsymbol{u}(t, \boldsymbol{\mu}) = \Phi (\boldsymbol{u}_0;t, \boldsymbol{\mu})$ with $\Phi: \mathcal{M} \times \mathbb{R} \times \Omega \rightarrow \mathcal{M}$ the solution map, $\mathcal{M}$ the solution manifold and $\Omega \subseteq \mathbb{R} ^ p$, the parameters space.
\begin{definition}[Lyapunov exponent]  
\label{def: 1}
Considering a matrix function $d_{\boldsymbol{u}}\Phi$ and a nonzero vector $\boldsymbol{w}$ of the Euclidean space $\mathbb{R}^n$ the quantity
\begin{equation}
\lambda(d_{\boldsymbol{u}}\Phi,\boldsymbol{w})=\limsup_{t\to\infty}\frac{1}{t}\ln\|d_{\boldsymbol{u}}\Phi \boldsymbol{w}\|
\end{equation}
is called the \emph{1-dimensional Lyapunov Characteristic Exponent} or the 
\emph{Lyapunov Characteristic Exponent of order 1} (1-LCE) of $d_{\boldsymbol{u}\Phi}$ with respect 
to vector $\boldsymbol{w}$.
\end{definition}


Suppose that we introduce a small perturbation (of magnitude $\varepsilon_0$), along a direction $\boldsymbol{t} \in T_{\boldsymbol{u}} \mathcal{M}$, that we call $\boldsymbol{v}_0 = \varepsilon_0 \boldsymbol{t}$. We can study how the perturbation evolves along the trajectory $\boldsymbol{u}(t; \boldsymbol{\mu})$ by linearizing the dynamical's system governing equation around the reference solution $\boldsymbol{u}(t; \boldsymbol{\mu})$. By  expanding in Taylor series the flow $\boldsymbol{f}(\boldsymbol{u}(t, \boldsymbol{\mu}))$ the perturbed dynamical system can be written as
\begin{equation*}
   \cancel{\dot{\boldsymbol{u}}(t, \boldsymbol{\mu})} + \dot{\boldsymbol{v}}(t) =  \cancel{\boldsymbol{f}(\boldsymbol{u}(t, \boldsymbol{\mu}))} + d_ {\boldsymbol{u}} \boldsymbol{f}(\boldsymbol{u}(t, \boldsymbol{\mu})) \boldsymbol{v} + o(\varepsilon \boldsymbol{1}),
\end{equation*}
and so 
\begin{equation}
    \label{eq: variational equation}
    \dot{\boldsymbol{v}}(t) \approx  d_ {\boldsymbol{u}} \boldsymbol{f}(\boldsymbol{u}(t, \boldsymbol{\mu})).
\end{equation}
Equation \ref{eq: variational equation} is called \emph{variational equation} and it is a linear differential equation in the variable $\boldsymbol{v}$ with time-dependent coefficient matrix $d_ {\boldsymbol{u}} \boldsymbol{f}(\boldsymbol{u}(t, \boldsymbol{\mu}))$. The solution to equation \ref{eq: variational equation} can be written using the solution map $\Phi$ as:
\begin{equation*}
    \boldsymbol{v}(t) = d_ {\boldsymbol{u}} \Phi(\boldsymbol{u}_0; t, \boldsymbol{\mu}) \boldsymbol{v}_0.
\end{equation*}

By tracking the evolution of one perturbation, you will get with probability one, the largest LE \cite{skokos2009lyapunov}. This is the most important LE, since its sign is indicative of the nature of the system under investigation. The Multiplicative Ergodic Theorem (MET) guarantees the existence of the LEs and of a base $\{\boldsymbol{e}_1, \boldsymbol{e}_2, ..., \boldsymbol{e}_n\}$ of $T_{\boldsymbol{u}_0}\mathcal{M}$, such that the LEs associated to each element of the base \cite{oseledec1968multiplicative}:
\begin{equation*}
    \lambda(d_{\boldsymbol{u}}\Phi,\boldsymbol{e}_1) \ge \lambda(d_{\boldsymbol{u}}\Phi,\boldsymbol{e}_2) \ge ... \ge \lambda(d_{\boldsymbol{u}}\Phi,\boldsymbol{e}_n)
\end{equation*}
The collection of all the LEs form the Lyapunov Spectrum. The Lyapnunov Spectrum provides useful information about the inertial manifold and about the global attractor embedded within it.

\begin{definition}[Global Attractor]
\label{def: Global Attractor}
Let $H$ be a Hilbert space and $\{\Phi(\boldsymbol{u}_0, t)\}_{t \ge 0}$ the solution map associated to a dynamical system on $H$.  
A set $\mathcal{A} \subset H$ is called a \emph{global attractor} if:

\begin{enumerate}
    \item $\mathcal{A}$ is compact in $H$;
    
    \item $\mathcal{A}$ is invariant under the flow, i.e.,
    \[
    \Phi(\mathcal{A}, t) = \mathcal{A}, \quad \forall t \ge 0;
    \]
    
    \item $\mathcal{A}$ attracts all bounded subsets of $H$, meaning that for every bounded set 
    $\mathcal{B} \subset H$,
    \[
    \operatorname{dist}(S(t)\mathcal{B}, \mathcal{A})
    := \sup_{\boldsymbol{u}_0 \in \mathcal{B}} \inf_{y \in \mathcal{A}} \|\Phi(\boldsymbol{u}_0, t) - y\|_H 
    \longrightarrow 0 \quad \text{as } t \to \infty.
    \]
\end{enumerate}

In particular, for every initial condition $\boldsymbol{u}_0 \in H$, the trajectory satisfies
\[
\operatorname{dist}(S(t)\boldsymbol{u}_0, \mathcal{A}) \to 0 \quad \text{as } t \to \infty.
\]

\end{definition}

\begin{definition}[Inertial Manifold]
\label{def: Inertial Manifold}
Let $H$ be a Hilbert space and $\{\Phi(\boldsymbol{u}_0, t\}_{t \ge 0}$ the solution map associated to a dynamical system on $H$.  
A finite-dimensional smooth manifold $\mathcal{M} \subset H$ is called an \emph{inertial manifold} if:

\begin{enumerate}
    \item $\mathcal{M}$ is positively invariant under the flow, i.e.,
    \[
    \Phi(\mathcal{M}, t) \subset \mathcal{M}, \quad \forall t \ge 0;
    \]
    
    \item $\mathcal{M}$ attracts all trajectories at an exponential rate, i.e., for every 
    $\boldsymbol{u}_0 \in H$, there exist constants $\kappa_1, \kappa_2 > 0$ such that
    \[
    \operatorname{dist}(\Phi(\boldsymbol{u}_0, t) , \mathcal{M}) 
    \le \kappa_1 e^{-\kappa_2 t}, \quad \forall t \ge 0.
    \]
\end{enumerate}

\end{definition}

While definition \ref{def: Global Attractor} and \ref{def: Inertial Manifold} looks similar they are quite different mathematical objects. Indeed, while the global attractor could be a fractal set, the inertial manifold is required to be smooth, or Lipschitz. Another important difference concerns the rate at which trajectories that do not initially lie on these manifolds approach them. For the global attractor, the distance between a trajectory and the attractor is known to vanish asymptotically as time tends to infinity. In contrast, trajectories approach the inertial manifold at an exponential rate in time.. Last but not least, the global attractor is a subset of the inertial manifold.

For some dissipative systems (KS equation, Ginsburg-Landau equation, Cahn-Hilliard equations), it can be proved the existence of a inertial manifold \cite{temam1990inertial}. The Lyapunov exponents are tightly related with these manifolds. For instance, the Kaplan-York (KY) dimension gives an upper bound to the fractal dimension of the global attractor.
\begin{equation}
    \label{eq: KY dimension}
    D_{KY} = k + \frac{\sum _ {i=1} ^k \lambda_i}{|\lambda _{k+1}|}
\end{equation}

Where $k$ is the index such that $\sum _ {i=1} ^k \lambda_i \ge 0$ and $\sum _ {i=1} ^{k + 1} \lambda_i < 0$. 

A chaoticity measure of the dynamical system is given by the Kolmogorov-Sinai entropy, which can be interpreted as the average amount of new information generated by the system per unit of time, which means how unpredictable is the evolution of the system itself \cite{frigg2004sense}. The Pessin's theorem indeed states that the Kolmogorov-Sinai entropy can be expressed as the sum of the positive LEs \cite{skokos2009lyapunov}.

\section{FOM LEs of KS equation}
\label{sec: FOM LEs}
The variational equation for the KS system represented by Equation \ref{eq:ks_equation} can be written as:
\begin{equation}
    \label{eq: variational equation KS}
    \frac{\partial  v}{\partial t} +  \frac{\partial u v} {\partial x} + \frac{\partial^2 v}{\partial x ^2} + \nu \frac{\partial^4 u}{\partial v ^4}= 0.
\end{equation}
Both equation \ref{eq:ks_equation} and equation \ref{eq: variational equation KS} are solved together, using a spectral scheme coupled with ETDRK4 time stepping scheme. Specifically, in order to evaluate the non linear term in between the time discretization nodes, as ETDRK4 requires, we integrated equation \ref{eq:ks_equation} by using half of the time steps used to solve \ref{eq: variational equation KS}. Then, the LEs for the FOM are computed using the algorithm described by \cite{skokos2009lyapunov}. A validation of this solver was performed by comparing the computed spectra with the reference one reported in \cite{ozalp2025stability}, for the same value of the domain length $L$. The validation is not reported here. 

\section{Equivalence of LE}
\label{sec: Equivalence LE}
To show the equivalence of the LEs found by the FOM and ROM system's governing equations we consider the same framework introduced in section \ref{sec: Lyapunov Exponents}.
Suppose that 
\begin{equation*}
    T_{\boldsymbol{u}} \mathcal{M} = \text{span}(\{\boldsymbol{\xi}_1(\boldsymbol{u}),\, ..., \, \boldsymbol{\xi}_l(\boldsymbol{u})\})
\end{equation*}
where $\boldsymbol{\xi}_i(\boldsymbol{u})$ is the $i$-th covariant Lyapunov vector associated with the $i$-th Lyapunov exponent. Let us denote $E:\mathcal{M} \rightarrow \Sigma \subseteq \mathbb{R}^l$ a (full rank) manifold parameterization, and $\varphi = E \circ \Phi \circ E^{-1}$, the solution operator in the latent space. It is possible to apply Definition \ref{def: 1} to the solution map defined by $\varphi$ in the latent space: 
\begin{equation}
\label{eq: reduced lambda}
\lambda ^r (d_{\boldsymbol{u}}\varphi,\boldsymbol{v})=\limsup_{t\to\infty}\frac{1}{t}\ln\|d_{\boldsymbol{z}}\varphi \boldsymbol{v}\|.
\end{equation}

Let $\boldsymbol{w} \in T_{\boldsymbol{u}} \mathcal{M}$, and hence $\boldsymbol{v} = d_{\boldsymbol{u}} E \boldsymbol{w}$. We want to prove that  $\lambda ^r (d_{\boldsymbol{u}}\varphi,\boldsymbol{v}) = \lambda  (d_{\boldsymbol{u}}\Phi,\boldsymbol{w})$. We can then expand the terms inside the norm in Equation (\ref{eq: reduced lambda}) as:
\begin{equation}
\label{eq: norm}
    \| dE d\Phi d E^{-1} dE \boldsymbol{w} \|,
\end{equation}
since $\boldsymbol{w} \in T_{\boldsymbol{u}} \mathcal{M}$, $\boldsymbol{w} = d E^{-1} dE \boldsymbol{w}$.  Let $\sigma_{min}(\boldsymbol{u})$ and $\sigma_{max}(\boldsymbol{u})$ be maximum and minimum modules of the singular values of the linear map $dE$ at the point $\boldsymbol{u} \in \mathcal{M}$. Let us introduce the definitions
\begin{equation*}
    \overline{\sigma}_{min} = \inf_{\boldsymbol{u}\in \mathcal{M}} \left( \sigma_{min}(\boldsymbol{u}) \right) \quad \mathrm{and} \quad 
    \overline{\sigma}_{max} = \sup_{\boldsymbol{u}\in \mathcal{M}} \left( \sigma_{max}(\boldsymbol{u}) \right).
\end{equation*}
Being $E$ full rank, it is guaranteed that $\overline{\sigma}_{min}>0$ and a linear map between finite spaces is always bounded. 
The argument of the logarithm in equation (\ref{eq: reduced lambda}) can then be bounded as
\begin{equation*}
 \overline{\sigma }_{min} \|d\Phi \boldsymbol{w} \|  \le \|dE d\Phi \boldsymbol{w} \| \leq \overline{\sigma }_ {max} \|d\Phi \boldsymbol{w} \|.
\end{equation*}
By applying the logarithm to the above expression and exploiting the monotonicity of the logarithm, the same inequalities hold:
\begin{equation*}
 \ln (\overline{\sigma} _ {min} \|d\Phi \boldsymbol{w} )\|  \le \ln( \|dEd\Phi \boldsymbol{w} \|) \leq \ln(\overline{\sigma} _ {max} \|d\Phi \boldsymbol{w} \|).
\end{equation*}
Since the logarithm of a product is the sum of the logarithm of each factor, by the squeezing theorem, we can prove that $\lambda ^r (d_{\boldsymbol{u}}\varphi,dE\boldsymbol{w}) = \lambda  (d_{\boldsymbol{u}}\Phi,\boldsymbol{w})$:
\begin{equation*}
 \limsup_{t\to\infty}\frac{1}{t} \ln (\overline{\sigma} _ {min} \|d\Phi \boldsymbol{w} \| ) \le \limsup_{t\to\infty}\frac{1}{t} \ln( \|d\varphi \boldsymbol{v} \|) \\ \leq  \limsup_{t\to\infty}\frac{1}{t} \ln(\overline{\sigma} _ {max} \|d\Phi \boldsymbol{w} \|)
\end{equation*}
\begin{equation*}
\lambda  (d_{\boldsymbol{u}}\Phi,\boldsymbol{w}) \le\lambda ^r (d_{\boldsymbol{u}}\varphi,dE\boldsymbol{w}) \leq  \lambda  (d_{\boldsymbol{u}}\Phi,\boldsymbol{w})
\end{equation*}
\section{Symmetries and Lyapunov Exponents in Infinite-Dimensional Systems}
\label{sec: Symmetries and LEs}
Let us consider an abstract evolution equation on a Banach (or Hilbert) space $\mathcal{H}$:
\begin{equation}
    \frac{\partial u}{\partial t} = F(u), \qquad u(t) \in \mathcal{H},
\end{equation}
where $F : \mathcal{H} \to \mathcal{H}$ is a sufficiently smooth (e.g.,\ $C^1$) operator.
We denote by $\Phi(t,u_0)$ the associated solution map, so that $u(t)=\Phi(t,u_0)$ solves the equation with initial condition $u(0)=u_0$.

Let $\{g_\alpha\}_{\alpha \in \mathbb{R}}$ be a one-parameter group of transformations
\begin{equation}
    g_\alpha : \mathcal{H} \to \mathcal{H},
\end{equation}
such that:
\begin{align}
    g_0 &= \mathrm{Id}, \\
    g_{\alpha+\beta} &= g_\alpha \circ g_\beta,
\end{align}
and the map $\alpha \mapsto g_\alpha u$ is $C^1$ for all $u \in \mathcal{H}$.

We say that the evolution equation is \emph{equivariant} under this group if
\begin{equation}
    F(g_\alpha u) = g_\alpha F(u), \qquad \forall u \in \mathcal{H}, \ \alpha \in \mathbb{R}.
\end{equation}

Equivalently, the flow satisfies
\begin{equation}
    \Phi(t, g_\alpha u_0) = g_\alpha \Phi(t, u_0),
    \qquad \forall t \in \mathbb{R}.
\end{equation}

Since the symmetry is continuous, we define the infinitesimal generator
\begin{equation}
    v(u) := \left.\frac{d}{d\alpha} g_\alpha u \right|_{\alpha=0}.
\end{equation}

Under the above assumptions, $v : \mathcal{H} \to \mathcal{H}$ is a well-defined vector field tangent to the group orbit $\{ g_\alpha u : \alpha \in \mathbb{R} \}$. Differentiating the equivariance relation $\Phi(t, g_\alpha u_0) = g_\alpha \Phi(t, u_0)$ with respect to $\alpha$ at $\alpha=0$, we obtain
\begin{equation}
    d_u \Phi(t,u_0)\, v(u_0)
    =
    v(\Phi(t,u_0)).
\end{equation}

Thus, the vector
\begin{equation}
    v(t) := v(u(t)) = v(\Phi(t,u_0))
\end{equation}
is a solution of the linearized equation
\begin{equation}
    \dot{v} = dF(u(t))\, v.
\end{equation}

In particular, the subspace spanned by the symmetry generator is invariant under the linearized dynamics.

The Lyapunov exponent associated with this direction is defined as
\begin{equation}
    \lambda_v(u_0)
    =
    \limsup_{t \to \infty}
    \frac{1}{t}
    \log \| D_u \Phi(t,u_0)\, v(u_0) \|.
\end{equation}

Using the identity above, this becomes
\begin{equation}
    \lambda_v(u_0)
    =
    \limsup_{t \to \infty}
    \frac{1}{t}
    \log \| v(\Phi(t,u_0)) \|.
\end{equation}

Therefore, the Lyapunov exponent associated with the symmetry direction is determined entirely by the growth of the norm of the generator along the trajectory.

In particular, if there exist constants $C > 0$ and $p \ge 0$ such that $\| v(\Phi(t,u_0)) \| \le C (1+t)^p$,
then, as shown in\cite{pikovsky2016lyapunov}, $\lambda_v(u_0) = 0$.

\section{NNs hyper parameters}
Hyperparameter tuning for the autoencoder was performed exclusively on the viscous Burgers’ equation, with the resulting architecture subsequently applied to both the linear advection and Kuramoto–Sivashinsky (KS) problems. Given the relative simplicity of the linear advection case, we set the function $f$ in Equation \ref{eq: projecion layer} to zero yielding a weight-free decoder. For the more complex viscous Burgers’ and KS cases, we defined $f$ as a ConvTranspose1d layer from PyTorch
\begin{table}[ht]
\centering
\label{tab:cae_hyperparameters}
\begin{tabular}{ll}
\hline
\textbf{Hyperparameter} & \textbf{Value} \\
\hline
Input size & $N_h$ \\
Latent dimension & $l$ \\
Kernel sizes & $[21, 21, 21, 9, 5]$ \\
Channels & $[1, 2, 4, 8, 8, 8]$ \\
Strides & $[4, 4, 2, 2, 2]$ \\
Padding mode & \texttt{circular} \\
Bias & \texttt{True} \\
\texttt{Number Linear Layers} & $1$ \\
\hline
\end{tabular}
\caption{Hyperparameters of the convolutional autoencoder. The activation function for the classical CAE is the Tanh activation function, while for the symmetric case we choose the Rational activation function (Equation \ref{eq: rational AF}) with $\alpha=0.5$.}
\end{table}

\begin{table}[ht]
    \centering
    \begin{tabular}{ll}
    \hline
    \textbf{Test case} & \textbf{Layers size}   \\
        \hline
         Linear Advection & $[2, 4, 8, 8, 4, 2]$ \\
         Viscous Burgers & $[3, 8, 16, 16, 8, 2]$ \\
         KS & $[24, 48, 96, 192, 192, 48, 24]$\\
         \hline
    \end{tabular}
    \caption{Architecture of the LF network for each test case. The LF network is a feed-forward neural network with Tanh activation functions. The input, hidden, and output layer sizes are reported in the table. }
    \label{tab:LF net table}
\end{table}

\end{document}